\documentclass[12pt,a4paper]{amsart}
\usepackage{amssymb,amsmath,amsthm}
\usepackage[utf8]{inputenc}
\usepackage{amsmath}
\usepackage{MnSymbol}
\usepackage{tikz}
\usepackage{tikz-qtree}
\usepackage{epigraph}
\usetikzlibrary{positioning,shapes}

\usepackage{hyperref}
\newtheorem{thm}{Theorem}

\newtheorem{definition}[thm]{Definition}
\newtheorem{srule}[thm]{Rule} 
\newtheorem{theorem}[thm]{Theorem}

\def\lall{\,  \{  }
\def\rall{  \} \, }
\def\lexi{\,  [ }
\def\rexi{ ]\, }
\def\prop{:}
\def\eq{=}

\def\Nf{N}
\def\0f{{0}}
\def\zf{{\bf 0}}
\def\1f{{\bf 1}}
\def\of{{\bf 1}}
\def\nf{{ s}}
\def\lef{\preceq}
\def\lf{ \prec }

\newcommand{\N}{\mathbb{N}}
\newcommand{\X}{\mathbb{X}}
\newcommand{\Y}{\mathbb{Y}}
\newcommand{\I}{\mathbb{I}}

\newcommand{\nonmem}[2]{#2 ( #1 ) = #2} 
\newcommand{\req}[1]{#1=#1}
\newcommand{\ds}[1]{#1'}
\newcommand{\mem}[2]{#2 ( #1 ) \ds{\eq} #2}

\begin{document}

\title{A formal system of mathematics based on definitions}
\author{Christoph Thiele}

\address{Mathematical Institute, University of Bonn, Endenicher Allee 60, 53115 Bonn, Germany}
\email{thiele@math.uni-bonn.de}


\date{\today}
\maketitle

\section{Introduction}
\label{section:introduction}

It has been understood for more than a century that mathematical reasoning can in principle be 
formalized to the extend that it becomes machine readable and verifiable. 
The implementation of this program has gained momentum in the age 
of efficient computers. Many theorems have been translated into 
machine readable code in recent decades. The challenges of 
the implementation raised the desire for code that is convenient both 
for computer and humans. For example,  the language ForTheL 
\cite{paskevich} uses English phrases ubiquitous in mathematics 
to assemble English sentences in its code.

The goal of this paper is  to assemble a formal system out of a minimal
collection of typical ingredients in a mathematics paper. Instead of English text, we 
write a stenographic notation that follows closely such ingredienrs but resembles more 
classical mathematical formulas. We then use this formal system for some first exploration in mathematics in Sections \ref{section:basic}  and \ref{section:peano}. We establish 
Russell's paradox, construct the natural numbers and prove the Peano axioms.

To compare our system with standard axiomatizations in mathematics, recall that
working mathematicians generally accept a standard of mathematical axiomatization
based on propositional and predicate logic. These types of logic need an underlying 
content. The most commonly accepted universal content is Zermelo-Fraenkel set theory.
A second example is Peano arithmetic.  Basic statements 
of set theory concern sets being elements of other set, while basic statements in Peano 
arithmetic concern algebraic equations. Our system introduces a form of propositional 
and predicate logic and links it from the beginning to a content about objects and
their equality. Further content arises in the process of doing mathematics. The identity 
of an object may not be entirely determined at the beginning of its use, but evolve 
with choices made throughout the mathematical process. This  allows to build a system 
sophisticated enough to match Zermelo-Fraenkel set theory . We also develop a similar 
logical structure around attributes that are constructed from the basic attribute of equality. 
This structure is more basic and more constructive than propositional and predicate logic 
but equally important to describe the the growing material developed in mathematics.
 
 In Section \ref{section:formal-rules}, we give a concise formal description of our system.
 As an informal introduction, we discuss here the elements of the system and 
 compare them with typical mathematical  concepts.
 Basic statements are built from an attribute and a number of objects.
 The link between attribute and objects is typically expressed by the words ``is'' or "are".
An example of a basic statement is
$$Nx$$
where the capitalized letter stands for an attribute and the lower case letter stands for an object.
For example, the letter $N$ may have an interpretation as the attribute of being a natural number, 
then the formula means that $x$ \emph{is} a natural number.
More precisely, there is a temporal component in our code. Reading the above as ``is'' assumes that
 the object $x$ has been introduced earlier in the argument. If $x$ does not occur earlier than the 
 above line,  then this line introducces the object $x$ and can be interpreted as saying \emph{let} $x$ \emph{be} a natural number. 
 
Attributes may also take more than one object, the fundamental example is
$$x\eq y$$
meaning that $x$ and $y$ \emph{are} {equal}, or $x$ is {equal} \emph{to} $y$.  
The latter formulation stresses more the order of the objects, but in case of the symmetric equality 
this emphasis is not necessary. The equality sign is the basic attribute
within our formal system. It has the structural effect that it allows one of the equal 
objects to be substituted by the other in any further statement, Rule \ref{rule:subf}.
The link between equality and substitution appears as early as in Freiherr von  Wolff's \cite[\S 17]{wolff}, ``Wenn ich ein Ding B f\"ur das
Ding A setzen kann und es bleibet alles wie vorhin, so ist
A und B einerlei'' and  ``Wenn ich aber B f\"ur A setze kann und es
bleibet nicht alles wie vorhin, so sind
A und B unterschieden oder verschiedene Dinge.''
A void substitution is caused by having the same variable on both sides of an equality.
Such reflexive eqalities are therefore true by default.

The negation ``not'' comes with its own symbol, a prime symbol behind an attribute.
We write
$$x\ds{\eq}y$$
for $x$ and $y$ are \emph{not} equal.
Indeed, all attributes come in pairs, one with and one without the
prime, which are classically dual in the sense of 
negation.
We adopt the principle of the excluded middle, so we allow
arguments by case distinction between a statement and its dual
statement, Rule \ref{rule:cases}. We also adopt the rule of explosion, that is,
a statement is proved by bringing its negation to a contradiction,
Rule \ref{rule:contradiction}. 

We represent the phrase ``there exists" by brackets as in
$$\lexi N\xi \rexi \xi=o$$
meaning that \emph{there exists} a natural number $\xi$ such that $\xi$  is equal to $o$.
Similarly, we represent the phrase ``for all''  as
$$\lall N\xi \rall \xi\neq o$$
meaning that  \emph{for all} natural numbers $\xi$ we have that $\xi$ is not
equal to $o$. Typographically, we have chosen brackets common on standard keyboards.
The existential brackets are reminiscent of the letter ``E''
and the universal brackets are reminiscent of the letter ``A''.

It appears that exactly one of the above bracketed statements is true. Indeed, these statements are
formally dual to each other.
Truth depends on what the object $o$ might be. If for example it is zero, then
the former statement is true and the latter is false, while if $o$ is not a natural number,
the latter  is true and the former is false.

The grammar explicitly suggests that $\xi$ is an object that we maintain some liberty to
more narrowly specify later.  We use greek letters for such deferred objects, and latin letters for objects that 
we deem determined already or that we determine at the time of first use.
Such distinction of variables occurs in Frege's Begriffsschrift \cite{frege}:
``Alle Zeichen, die ich
verwende, theile ich daher in
solche, unter denen man sich
Verschiedenes vorstellen
kann, und in solche, die einen
ganz bestimmten Sinn haben.''

Replacing greek by latin letters, we have similar statements to the above. The statement
$$\lexi Nx \rexi x\eq o$$
means that $x$ is a  natural number \emph{and}  $x$ is equal to $o$, while
$$\lall Nx \rall x \ds{\eq}o$$
means that \emph{if} $x$ is a natural number, \emph{then} $x$ is not equal to $o$. 
The type of letter switches therefore between the quantifier phrases``there exists'' 
and ``for all'' on the one hand and the logical connectives ``and'' and ``if, then'' on the other hand.

These quantifier brackets are formally a process of concatenating statements.
Rules \ref{rule:exist} and \ref{rule:all} regulate this concatenation process for these
brackets. In the  case of existential brackets, we can simply deduce 
both parts of the concatenation from the concatenated statement. If a greek letter occurs
at suitable place we use the moment of splitting the existential statement to introduce an object  
substituting a greek by a latin letter. In the case of universal brackets, we need that the 
first of the concatenated statements appears in the past, again with a 
latin letter in place of a greek letter when necessary, and we  are then allowed to
deduce the second statement, possibly with the same modification of the letters.

The deferment of determination of an object leads very naturally to dependence 
of objects of each other. Consider the statement
$$\lall N\xi \rall \lexi N\eta \rexi \eta \ds{\eq} \xi$$
meaning that for all natural numbers there is another natural number  that
is not equal to the former. This is a statement that appears true. For example,
if the first number is zero, we pick the second number  to be one, and if the first
number is not zero, we pick the second number to be zero. Note the importance of 
the order of determining the two numbers.
We need to know what the first number is, before we can safely specify the second. 
It is therefore problematic, to replace in the above formula only the second greek letter by a latin letter,
$$\lall N\xi \rall \lexi Ny \rexi y \ds{\eq} \xi$$
This assumes we have fixed the latin object  already and it means that for
all natural numbers  we have that the fixed object is a natural number and the natural number
 is not equal to the fixed object. If the first conclusion, namely that the fixed object is a natural number, is correct, then the second conclusion is incorrect if we choose the natural number to be the same as 
 the fixed object. We adopt instead the dependence notation 
$$\lall N\xi \rall \lexi Ny(\xi) \rexi y(\xi) \ds{\eq} \xi$$
where $y(\xi)$ may be read as $y$ \emph{of} $\xi$.
So our system has a functional bracket for concatenating objects into new objects.
Traditionally, the first object in the concatenation is interpreted as a function that associates to
each possible argument, that is the second object in the concatenation, a new object.
The function may well be determined prior to the determination of the argument.
The function in the last displayed statement is called a choice for the  earlier displayed statement
with corresponding greek variable, which is a guarantor for the possibility of a suitable choice. The passage from 
the statement with a greek letter to the statement with a concatenated object
is close to the classical axiom of choice in Zermelo-Fraenkel set theory. 
In our system, choice lies at the heart of the formal system, Theorem \ref{thm:choice}.
However, the replacement by a concatenated object has to be done with some care. 
Consider the true statement
$$\lall  \xi= \xi \rall \lexi  \eta= \eta \rexi 
\lexi  \lall  \xi(\xi)\eq \xi\rall \eta\ds{\eq} \xi \rexi 
\lall  \xi(\xi)\ds{\eq} \xi\rall \eta {\eq}  \xi$$
which states that for every object there
is another object  that is either equal or not equal to  first object depending on
a certain condition on the first object. 
Substituting the second greek letter by a concatenated object, we obtain
$$\lall  \xi= \xi \rall \lexi  y(\xi)= y(\xi)\rexi 
\lexi  \lall  \xi(\xi)\eq \xi\rall y(\xi)\ds{\eq} \xi \rexi 
\lall  \xi(\xi)\ds{\eq} \xi\rall y(\xi) {\eq}  \xi$$
This turns out to be a false statement. It makes a claim for an arbitrary object, but the last part
of the statement fails  if we let this object be the newly introduced object $y$. Namely, we obtain
$$
\lexi  \lall  y(y)\eq y\rall y(y)\ds{\eq} y \rexi 
\lall  y(y)\ds{\eq} y\rall y(y) {\eq}  y$$
which states that certain two objects should be equal precisely if they are not equal.
The problem is that the  statement prior to applying choice works for all objects known 
prior to applying the choice. The  choice then creates a new object, and this new object 
refuses to abide by the universal statement used to create itself. The example is a close variant of
Russell's paradox in \cite{russell}, which revealed a fault in this context in Frege's 
Begriffsschrift \cite{frege}. The example is also somewhat in the spirit of the colloquial paradox 
"This sentence is a lie". 
Theorem \ref{thm:choice} describes a more careful choice, avoiding that the new object 
will be forced to satisfy the statement used to create itself.

Within our system, we interpret objects and  concatenated object at the same level
and are in particular free to equate and potentially substitute any such objects.
Statements that play a particular role in our system are the statements
$$y(x)\eq y$$
$$y(x)\ds{\eq} y$$
They appear somewhat odd in classical mathematics due to the type mismatch
between function and value of the function. 
In our system. The statements encode what is classically considered the domain of a function.
The second statement says that the argument of the concatenation is in the domain of the 
function and the first says that the argument is not in the domain of the function.
Note how typographically it appears in the first statement that the function
refuses to accept the argument and just remains itself.

The concept of domain is used in the formalism of choice  to avoid Russell's paradox. 
Validity of universal statements  is only required for objects in the domain of some other
object. To have an interesting theory, one then needs to guarantee objects with somewhat 
large domains. This is established by Theorems \ref{thm:universe} and \ref{thm:infinity}
which corresponds to the classical Zermelo-Fraenkel axioms
of replacement, union, and power set on the one hand and the axiom of infinity on the other
hand.

The deferment of the identity  of an object until certain other objects are identified 
allows to build a rich structure of 
interdependent objects. Functions are used in lieu of sets at the center 
of attention by von Neumann \cite{vNeumann}, who justifies 
his choice by practicality:
``Die technische Durchf\"uhrung
gestaltet sich jedoch beim
Zugrundelegen des Funktionsbegriffes wesentlich
einfacher, allein aus diesem
Grunde haben wir uns f\"ur
denselben entschieden.''
The rich structure of functions matches the richness of the classical von Neumann
universe in set theory. For convenience, and in analogy with the axiom of extension 
in the Zermelo-Fraenkel set theory, Theorem \ref{thm:function-uniqueness} explicitly restricts 
attention to how objects react with the choice of other objects. Objects  do not have additional characteristics of their own.

We turn our attention to attributes, which are represented by the equality symbol or by
capital letters. Unlike lower case letters, behind which there is a universe of objects
that  can be disacovered with the full arsenal of non-constructive mathematical 
reasoning such as proof by cases, by contradiction, and by choice, the
capital letters represent very concrete and explicitly constructed attributes. One should think
of attributes as abbreviations for more complicated statements build from the basic
attribute of equality. 

To introduce a new attribute, one has to
precisely write the abbreviated statement following Rule \ref{rule:con-pro}. 
We introduce for example the attribute of being a natural number by the string
$$N\xi\prop \lall S\phi\rall \mem \xi \phi $$
Here the statement to the left of the construction symbol $\prop$, which may represent the
phrase ``define'', is an abbreviation for the
statement to the right of the symbol. The right hand side uses more basic 
attributes such as the equality sign and capital letters that are defined earlier.
 The construction symbol is graphically a shorter variant of an equality symbol. Indeed, it  plays a 
 similar role in that it allows to substitute  the statement on one side of it by the statement on the other side  of it  in any further statement. Rule \ref {rule:subs} governs this process, which is subject
to some minor syntactical rules concerning brackets. 
The construction symbol is used solely as part of definitions, 
 in particular it does not come with a negated form of itself.

Despite being explicit, the attribute letters also come as latin and greek letters with similar
effect on deferment.
Latin attributes are assumed fixed already at the time of their use, while greek letters
may be used before they are replaced by some precisely defined latin attribute.
Greek attribute letters appear in documents called definitions and in claims of theorems, 
that are recorded for use in later documents. Definitions
are governed by Rule \ref{rule:des-def} defining objects based on earlier existence proofs
and Rule \ref{rule:con-def} defining attributes by construction. Claims are governed
by Rule \ref{rule:theorem}, these claims have to be followed by proofs except
in the small number of exceptional cases of Theorems \ref{thm:function-uniqueness},
\ref{thm:choice}, \ref{thm:universe}, and \ref{thm:infinity}, which have the
status of rules or axioms of the system.
Objects or attributes recorded in definition documents may depend on greek
attributes which will be specified later. This leads to a functional notation with
attributes as arguments. For example,
$$m(\Gamma)$$
could express the minimum of all natural numbers which satisfy the attribute $\Gamma$, and
$$A(\Gamma)\xi$$
could denote the attribute of $\xi$ being a natural number  larger than all natural natural numbers satisfying 
attribute $\Gamma$. The functional notation is purely a way of organizing deferred explicit 
constructions, there is no rule similar to a choice axiom leading to functional expressions with
attributes as arguments.
 
Unlike in the case of functional dependence of objects, where we have the liberty to change
the number of arguments by concatenating further functional brackets, there is a rigidity around attributes in terms of number of arguments owed to the explicitness of the construction. The number of arguments in the context of attributes is regulated by arities associated with objects or attributes.

\section{Acknowledgement}
The author acknowledges the contribution of many colleagues and students
in the development of this article.
The author acknowledges the support from the University of
California at Los Angeles and the University of Bonn,
as well as the Hausdorff Center for Mathematics.

\section{Formal rules}
\label{section:formal-rules}

\subsection{Documents}

Our  system of mathematics evolves as a growing number of
\emph{documents} of two {types}, \emph{definitions} and
\emph{theorems}. We \emph{label} the documents with natural numbers 
greater than nine, with larger numbers referring to later documents.  
The numbers zero to nine enumerate rules for the formal system.

Each document has a \emph{preamble}, containing the
{type} and label of the document,  an informal
{\it title} for mnemonic purpose but without direct relevance
to the formal rules, and finally a possibly empty list of further 
labels of selected earlier definitions.

The main body of the document is a collection of \emph{strings}.
A definition  contains a number of strings linearly ordered in time.
We write these strings into a chain of rectangular boxes from top
to bottom. 
A theorem consists of one earliest string called the \emph{claim},
and in most cases a \emph{proof}. The proof is a collection of  strings put into
rectangular boxes and arranged in a tree structure. Some
boxes may be empty. Two boxes are connected by at most one edge.
Each edge connects an earlier box with a later box called parent and child of each other.
We draw the parent higher than the child.
There is an earliest box of the proof called the \emph{root}.
Each box other than the root has exactly one parent.
Each box has zero, one, or two children. A box without children is
called a \emph{leaf}.
A box in a proof is earlier than another box, if it can be reached
from the latter by a path with each step going from a
box to its parent.
 
An \emph{ancestor} of a string in a given document is the string itself, or a string in the
same document that is earlier than the given string, or a string in a definition that is earlier than the given
document, or a claim in a theorem that is earlier than the given document.

Next to an edge we may write a justification for the edge. To the left of it we may write the number of one 
or several rules of the system, or labels of prior definitions  or theorems containing relevant ancestors. To the right
of an edge in a proof we may write relevant ancestors within the same proof, where each such ancestor
is identified by the number of steps from box to parent needed to reach the ancestor from the leaf. 
 
\subsection{Symbols}

Strings are finite chains of \emph{symbols}. 
Some possible symbols are the opening and closing \emph{functional brackets},
$$ ( \ \ \ )$$
the opening and closing  \emph{existential quantifier}  or \emph{existential} brackets,
$$ \lexi \ \  \ \rexi $$
the opening and closing  \emph{universal quantifier}  or \emph{universal} brackets,
$$ \lall \ \  \ \rall $$
the \emph{duality symbol}, and the \emph{construction symbol}.
$$ \ds{ }\ \ \  \prop $$

All other symbols are called letters.
There are \emph{object letters} and \emph{attribute letters}. Each of these come in two forms,
\emph{definite letters} and \emph{indefinite letters}.

We  use characters of the greek alphabet for indefinite letters and possibly add a natural number as subscript. 
Object letters are lowercase and attribute letters are uppercase. Examples of indefinite letters are 
$$ \xi \ \ \   \eta \ \ \ \xi_0 \ \ \ \xi_1\ \ \  \Gamma\ \ \  \Delta \ \ \ \Gamma_0 \ \ \ \Gamma_1$$
One particular definite attribute letter is the \emph{equality sign},
$$\eq $$
We use the latin alphabet for all other  definite letters and add a subscript containing the label of a document.
Object letters are lowercase and attribute letters are uppercase. Examples of definite letters are
$$x_{\ref{def:ref-eq}}\  \ \ x_{\ref{thm:function-uniqueness}}  \ \ \  y_{\ref{thm:function-uniqueness}}\ \ \ 
A_{\ref{def:ref-eq}}\ \  \ A_{\ref{thm:function-uniqueness}}  \ \ \  {B}_{\ref{thm:function-uniqueness}}
$$
A new combination of latin character and subscript in some document has the label of the present document as subscript.
Subscripts of latin characters therefore refer to the earliest document in which the combination of 
character and subscript occurs.
We may omit the subscript of a latin character inside a document, if the subscript is the earliest label among the 
present document's label and the labels listed in the preamble of the present document  such that this combination
of character and subscript appears in the document with that label. This is the purpose of the list of labels in the preamble.

With each definite letter, we associate a natural number called the \emph{attribute arity} of the letter.
With each definite attribute letter we additionally associate a pair of natural numbers called the
\emph{object arity} of the letter.
If the attribute arity of a definite letter is not zero, we further associate
with it for each number from one to the attribute arity a pair of numbers called the \emph{implicit arity} 
of the letter and the number.
The attribute arity of the equality sign is zero, its object arity
is one and one, and it has no implicit arity.

Given a string in a document, we associate to each indefinite attribute letter
occurring in the string an object arity consisting of two natural numbers.
The object arity of an indefinite attribute letter may change between different 
strings in a document.

We say that 
a letter is \emph{activated} at a string of a document, if it occurs in the string but not 
in any other ancestor of the string

\subsection{Terms}

We first concatenate symbols into  \emph{attributes} and \emph{terms}.

An attribute is an attribute letter or a concatenation
from left to right of an attribute letter and the duality symbol.
The two possible attributes containing a particular
attribute letter are called dual to each other.
We call an attribute definite or indefinite following the type of its attribute letter. 
We associate with an attribute the same arities as with the attribute letter in it.

A chain of consecutive symbols inside a string we call a substring, 
if it is not followed to the right by a duality symbol inside the larger string.
As duality symbols are only used directly following an attribute letter,
each attribute letter inside a string is part of a unique substring
which is an attribute.
  
An indefinite attribute or a definite attribute  with zero attribute 
arity is an \emph{attribute term}. 
Given a definite attribute of non-zero attribute arity,
we obtain further attribute terms by concatenating from left to right
first the given attribute, then the opening functional
bracket, then in succession as many attribute terms as the attribute arity of the
given attribute, and then the closing functional bracket. We require
that the object arity of each of these attribute terms is equal to
the implicit arity of the given attribute, which is the left most attribute
in the concatenation, and the number of the attribute term in the succession. 
The object arity of the concatenated attribute term is the object arity of its left most attribute. Two concatenated attribute terms are dual, if they are obtained form each other by replacing the left most attribute by its dual. 

A \emph{basic object term}  is an indefinite object letter, a definite object 
letter with zero attribute arity, or the concatenation from left to right
of a definite object letter with non-zero attribute arity, then
the opening functional bracket, then in succession as many attribute terms 
as the attribute arity of the object letter, and then the closing functional bracket. 
We require that the object arity of each of these attribute terms is equal to
the implicit arity of the object letter and the number of the attribute term in the succession.

A basic object term is an \emph{object term}. Given two not necessarily different
object terms, we may obtain a further object term by concatenating from left to 
right first one object term, then the opening functional bracket, then the other object term,
and finally the closing functional bracket.

Given any of the above terms, one may detect the most resent concatenation step of the
term.  The explicit brackets of the concatenation step are the right most closing bracket and the 
right most opening bracket that has the same number of opening and closing brackets between 
itself and the final closing bracket. If the most recent concatenation used a succession of
several attribute terms,  one may detect them as those attribute terms between the explicit concatenation brackets that do not lie between any 
further pair of opening and closing  brackets.

 \subsection{Statements}

Given an attribute term, we may concatenate an \emph{unquantified statement}
by concatenating in succession from left to right first as many object terms as the first number
of the object arity of the attribute term, then the attribute term, and then in succession
as many object terms as the second number of the object arity of the attribute term.
Given an unquantified statement, one may detect the components of this concatenation,
they are all the object and attribute terms which do not lie between any functional brackets
and are not followed to the right by an opening functional bracket.

Two unquantified statements are dual to each other, if they are obtained from each other
by replacing the attribute term used in the concatenation by its dual.
An unquantified statement is {\it admissible}, if it contains no indefinite object letters.

We call an unquantified statement \emph{equality}, if the attribute used in its concatenation
is the equality symbol.  We call the dual of 
an equality an \emph{inequality}. We call the equality or inequality \emph{reflexive}, if the two object terms used
in the concatenation are the same. 

Given two statements, we may concatenate a further statement by
writing first an opening quantifier bracket, then one statement,
then the matching closing quantifier bracket, and then the other
statement. We call such a concatenated statement quantified.
One can reconstruct the constituents of the most recent concatenation
by identifying the closing quantifier bracket in the concatenation
as the left most closing quantifier bracket that has the same number of opening and closing
quantifier brackets between itself and the initial opening bracket.
We call the first constituent term the \emph{hypothesis} and the
second constituent term the \emph{conclusion} of the quantified statement.
We call a quantified statement \emph{existential}, if the matching pair
of brackets used in the most recent concatenation are existential, 
otherwise we call it \emph{universal}.

Two quantified statements are dual to each other, if precisely one is universal 
and they have the same hypothesis but  mutually dual conclusions.
A quantified statement is admissible if either both hypothesis 
and conclusion are admissible or the hypothesis is unquantified 
and contains precisely one indefinite object letter, possibly at several positions,
and the conclusion becomes admissible if this indefinite letter
is replaced by a definite object letter with zero  attribute arity
in all positions where the former occurs.

We say an unquantified statement is contained in a statement, if it
it is a substring inside the statement that is not followed to the right
by an opening functional bracket.

Consider  for each indefinite attribute in a statement an attribute term  
with the same object arity as the indefinite attribute, so that mutually dual indefinite attributes 
correspond to mutually dual attribute terms and the attribute terms contain only definite letters.
We obtain a \emph{subordinate} of the 
statement by replacing each of the indefinite attributes in the statement by the corresponding 
attribute term in all positions where the former occurs.
Given a subordinate of a statement, and two indefinite object letters,
one occurring in the subordinate and the other not occurring in the subordinate,
we may obtain a further subordinate of the statement by replacing
the former indefinite letter by the latter in all positions where the former occurs.

\subsection{Constructions}

A \emph{construction} is concatenated from two statements by writing
from left to right one statement, then the construction symbol, and then the other statement.

We require that the left statement is unquantified
and it contains precisely one definite letter, namely the left most
letter in the attribute term used  in the concatenation of the unquantified statement.
We require that the left statement does not contain the duality symbol, 
and no letter occurs in more than one position in the left statement.
We require that all indefinite letters in the left statement also
occur in the right statement, but the definite letter in the left statement 
does not occur in the right statement.
We require that the right statement does not contain any indefinite attribute letters
other than those that also occur in the left statement.

We call the construction admissible if the right statement
becomes admissible when all indefinite object letters that appear in the left
statement are replaced by  a definite object letter of zero
attribute arity in all positions of the right statement where they occur.

Consider an admissible construction.
Consider  for each indefinite attribute in the construction an attribute term  
with the same object arity as the indefinite attribute, so that mutually dual indefinite attributes 
correspond to mutually dual attribute terms and the attribute terms contain only definite letters.
Consider for each indefinite object letter in the left statement an object term containing no indefinite
attributes. If the right statement is quantified, we also require the object term to not contain
any indefinite object letters.
We obtain a \emph{subordinate} of the construction by replacing each of the indefinite attributes 
in the statement by the corresponding 
attribute term in all positions where the former occurs and replacing each indefinite
letter occurring in the left statement by the corresponding object term in all positions 
of the construction where the former occurs. 
 
 Given a subordinate of a construction, we can obtain a further subordinate 
by replacing both left and right statements  by the respective dual statements.

Given a subordinate of a construction, and two indefinite object letters,
one occurring in the subordinate and the other not occurring in the subordinate,
we may obtain a further subordinate of the construction by replacing
the former indefinite letter by the latter in all positions where the former occurs.

A subordinate still has one construction symbol and statements on both sides
of the construction symbol. We call these statements  the two sides 
of the subordinate.

\subsection{Rules of adding strings}

The process of formal mathematics is to add strings to the most recent document, 
until it is complete, and to start a new document with a preamble when the previous
document is complete.

Assume the most recent document is a definition.
We can decide it is complete, or add a box with
a string at the bottom of the chain by 
one of the following two rules.

\begin{srule}[Construction in definition]
\label{rule:con-def}

We may add an admissible construction which activates the
definite letter in the left statement.

\end{srule}

\begin{srule}[Designation in definition]
\label{rule:des-def}

We may add a modified copy of the hypothesis of an existential claim of a prior theorem,
assuming this hypothesis is unquantified and contains a unique  indefinite letter.
Here modification means that this indefinite letter is replaced by 
one and the same object term in all positions where the former occurs. 

We require that this  object term has only one object letter,
only indefinite attribute letters, and no duality symbol.
The object letter is definite, activated at the new string,
and its attribute arity is equal to the number of different
indefinite attribute letters in the claim.
Each indefinite attribute letter in the claim appears in the object term
and has the same object arity in the new string as in the claim.
 
We write the label of the theorem at the new edge.

\end{srule}

Now assume the most recent document is a theorem.
If it has no claim yet, we proceed with the following rule.

\begin{srule}[Claim and root]
\label{rule:theorem}

We add a claim that is an admissible statement and
does not activate any definite letter.

If the theorem is one of Theorems \ref{thm:function-uniqueness} or
\ref {thm:choice} or \ref{thm:universe} or \ref{thm:infinity}
below, the theorem is complete with the claim.
Otherwise, we start a proof with a root that contains a subordinate of the dual of the claim,
where we require that each indefinite attribute letter in the claim
is replaced by a definite attribute letter that is activated at the root
and has the same object arity as the indefinite letter.
We require that identical indefinite attribute letters are replaced by
identical  letters and different indefinite attribute letters are replaced by
different letters.
\end{srule}

Now assume the most recent document is a theorem that already has a claim
and a root.  If all leaves of the proof are empty, the theorem is complete. If 
the proof has leaves that are not empty, we consider one such leaf 
and add one or several children to this leaf using one of the next seven rules.

\begin{srule}[Cases]
\label{rule:cases}

We add a pair of children to this leaf with admissible statements
that are mutually dual, contain no indefinite attribute letter, and do not
activate any definite letter. 

\end{srule}

\begin{srule}[Contradiction]
\label{rule:contradiction}

If a subordinate of an ancestors of the leaf is dual to the statement in the leaf,
we may add an empty child to this leaf and write the label of the ancestor at the new edge.

If the leaf has is a reflexive inequality, we  may add an empty child to the leaf. 

We may add a child to the leaf with an admissible
reflexive equality, which does not activate any  
attribute letter and activates at most one object letter.
We write the number of the present rule at the new edge.

If we intend to immediately after adding this child
use the reflexive equality for adding a child to this child using some other rule,
we may omit adding the first child and proceed to the next
child directly and write the present rule as part of the
justification for the child at the edge.

\end{srule}

\begin{srule}[Deduction]
\label{rule:all}

Consider a subordinate for  each of two ancestors of the leaf.
If one subordinate is universal and the other subordinate is the
hypthesis of the former, we may add
a child to the leaf with the conclusion of the former.

If one subordinate is universal with unquantified hypothesis
containing an indefinite object letter, and the other subordinate is the
modified hypthesis of the former, we may add
a child to the leaf with the modified conclusion of the former.
Here  modification means the said indefinite object letter is replaced by one and the same
object term in all positions where it occurs. 

Instead of adding the child as above, we may combine with a further application
of this Rule or Rule  \ref{rule:exist} to the quantified statement of the child, and
directly add a child according to the further application.

We write this rule and the relevant ancestors at the first new edge of the 
succession of applications of Rules \ref{rule:all} and  \ref{rule:exist}.

\end{srule}

\begin{srule}[Designation in proof]
\label{rule:exist}

Consider an object term containing no indefinite letters and a subordinate  of an ancestor of the leaf that is an existential statement.
Consider the modification of the subordinate obtained by doing nothing if the hypothesis of the subordinate is admissible and otherwise replacing the unique indefinite object letter occurring in the hypothesis  by the considered object term in all positions of the subordinate where the former occurs. 

We may add a child to the leaf with the hypothesis of the modified subordinate.
If we have done this, we may add a child to this child with the  conclusion of the modified subordinate.

We may omit the first child and directly add a child with the conclusion of the instantiation
to the original leaf, provided the hypothesis of the subordinate is admissible.

Instead of adding a child with the conclusion as above, we may combine with a further application of this Rule or Rule  \ref{rule:all} to the quantified conclusion,  and directly
add a  child according to the further application.

We write this rule and the relevant ancestor at the first new edge of the succession of applications of Rules \ref{rule:all} and  \ref{rule:exist}.

\end{srule}

\begin{srule}[Construction in proof]
\label{rule:con-pro}

We may add a child with an admissible construction that activates the
definite letter in the left statement and does not contain 
any indefinite attributes.

We write this rule at the new edge. 
\end{srule}

 \begin{srule}[Statement substitution]
\label{rule:subs}

Consider a subordinate of an ancestor that is a construction, and consider a statement, 
itself a subordinate of some ancestor of the leaf.

Assume one side of the subordinate appears as the right most substring of the 
statement. We may add  a child with a modification of the statement, meaning this rightmost
substring is replaced by the other side of the subordinate. 

Assume both sides of the subordinate are unquantified, 
and assume the statement  contains one side of the 
subordinate at one or several positions.
We may add a child  with a modification of the statement,
where modification means that this side of the  subordinate
is replaced by the other side of the subordinate in one or several positions 
where it occurs in the statement.

Instead of adding a child as above, we may combine with a subsequent applications
of Rule \ref{rule:exist} or Rule  \ref{rule:all}, breaking up the quantified statement 
of the child,  and directly add a child following the further applications.

We refer to this rule and the relevant ancestors in the first new edge of the 

\end{srule}

 \begin{srule}[Object substitution]
\label{rule:subf}

Consider a statement of the leaf 
and an indefinite object letter.

Consider two object terms which are the two object terms in the concatenation
of an equality that is an ancestor  of the leaf.
 
Assume a modification of the statement, where the indefinite letter
is replaced by one of the object terms in all positions, where it occurs,
is a subordinate of an ancestor of the leaf. We may add a child with 
another modification of the statement, where the indefinite letter is replaced 
by the other object term in all positions where the former occurs.

We write this rule and the relevant ancestors at  the new edge.
The modification may be void if the equality is reflexive, in which case we
copy the subordinate to the child and only write the relevant ancestor at the edge.

\end{srule}

\subsection{First definitions and theorems}

We begin the process of producing documents  until
the stage that the four theorems explicitly 
addressed in Rule \ref{rule:theorem} appear.

  \begin{theorem}[Object uniqueness]\ \\ \label{thm:function-uniqueness}
    $$\lall \phi= \phi \rall \lall \phi\neq \psi \rall  
    \lexi \phi(\xi)\neq \psi(\xi) \rexi \lall \nonmem \xi \phi \rall \mem \xi \psi  
    $$
   \end{theorem}

\begin{definition} [Domain by attribute] \ \\ \label{def:dom-prop}
\end{definition}

\begin{center}
\begin{tikzpicture}[every node/.style = {draw, shape=rectangle,anchor=north},growth parent anchor=south,level distance=.1 cm,sibling distance = 4.7cm]]
\node{$M(\Gamma) \prop \lexi \phi=\phi\rexi \lall \Gamma \xi \rall \mem \xi \phi $}
                    child{node{$D(\Gamma)\phi\prop \lexi \lall \Gamma \xi \rall \mem \xi \phi \rexi \lall \mem \xi \phi\rall \Gamma \xi $}
};
    \end{tikzpicture}
\end{center}

\begin{theorem}[Choice, \ref{def:dom-prop}]\ \\ \label{thm:choice}
    $$\lall M(\Gamma)\rall
    \lall \lall \Gamma \xi \rall \lexi \Delta \eta \rexi \eta \Sigma  \xi\rall
    \lexi     D(\Gamma) \eta \rexi
    \lall \Gamma \xi \rall \lexi \Delta \eta(\xi) \rexi \eta(\xi) \Sigma \xi
    $$
\end{theorem}

\begin{definition} [Range, union, power] \ \\ \label{def:range-union-power}
\end{definition}
\begin{center}
\begin{tikzpicture}[every node/.style = {draw, shape=rectangle,anchor=north},growth parent anchor=south,level distance=.1 cm,sibling distance = 4.7cm]]
\node{$R\phi \psi  \prop \lall \mem  \xi \phi \rall   \mem {\phi(\xi)} \psi$}
child{node{$U\phi \psi \prop \lall  \mem \xi \phi   \rall  \lall \mem \eta \xi \rall \mem \eta \psi$}
  child{node{$P\phi \psi \prop  \lall \nonmem \xi \psi\rall \lexi  \mem  \eta \xi \rexi \lall \mem \eta \phi \rall \xi(\eta)\ds{\eq} \eta $}
}
};
    \end{tikzpicture}
\end{center}

  \begin{theorem}[Universe, \ref{def:range-union-power}]\ \\ \label{thm:universe}
    $$ \lall \phi=\phi \rall  \lexi \mem \phi \psi\rexi \lexi R\phi \psi\rexi\lexi U\phi \psi\rexi  P\phi \psi
    $$
  \end{theorem}
 
\begin{definition} [Identity object, extension] \ \\ \label{def:id-ext-succ}
\end{definition}
\begin{center}
\begin{tikzpicture}[every node/.style = {draw, shape=rectangle,anchor=north},growth parent anchor=south,level distance=.1 cm,sibling distance = 4.7cm]]
  \node{$C\phi \prop \lall \mem \xi \phi \rall \phi(\xi)=\xi$}
  child{[level distance=.1cm]node{$E\phi \xi\prop  \lall \nonmem \xi \phi \rall \phi=\xi$}
    child{[level distance=.1cm]node{$H\psi \phi \prop
\lexi
\lall E\phi \xi\rall \mem \xi \psi
\rexi
\lexi
\lall \mem \xi \psi\rall E\phi \xi
\rexi C\psi$}
      }};
    \end{tikzpicture}
\end{center}

\begin{theorem}[Infinity, \ref{def:id-ext-succ}]\ \\
  \label{thm:infinity}
  $$ \lall \phi=\phi \rall  \lexi \mem \phi \psi\rexi
  \lall \mem \xi \psi\rall \lexi \mem \eta \psi \rexi H\eta \xi
  $$
  \end{theorem}

\section{Basic theorems and definitions}
\label{section:basic}

Equality is a binary attribute. In case of a reflexive equality, one may
want to replace this by a unary attribute. This is done in the following
definition.
\begin{definition} [Reflexive equality] \ \\ \label{def:ref-eq}
\end{definition}
\begin{center}
\begin{tikzpicture}[every node/.style = {draw, shape=rectangle,anchor=north},growth parent anchor=south,level distance=.1 cm,sibling distance = 4.7cm]]
\node{$I\xi\prop \xi=\xi$};
    \end{tikzpicture}
\end{center}

A number of rules of the system can be to similar effect be expressed
in the language of the system. We present some examples.
We first express that there is some object. It satisfies the attribute in Definition \ref{def:ref-eq}.
As a consequence of Theorem \ref{thm:existence-object}, a definite object letter may be activated
in a proof with this unary attribute, analoguous to Rule \ref{rule:contradiction}
for a reflexive equality, which is used in the proof of the theorem.

\begin{theorem}[Existence of object, \ref{def:ref-eq}]\ \\ \label{thm:existence-object}
  $$\lexi I\xi \rexi I\xi$$
\end{theorem}

  \begin{center}
\begin{tikzpicture}[every node/.style = {draw, shape=rectangle,anchor=north},growth parent anchor=south,level distance=.4 cm,sibling distance = 4.7cm]]
\node{$  \lall I\xi \rall \ds{I}\xi $}
                   child{node{$ x{\eq} x $}
                    child{node{$ I x $}
                    child{node{$ \ds{I} x $}
          child{node{} 
edge from parent node[left, draw=none] {\color{red}  \ref*{rule:contradiction}}
edge from parent node[right, draw=none] {\color{red} 1}
      }
                     edge from parent node[left, draw=none] {\color{red} \ref*{rule:all}}
      edge from parent node[right, draw=none] {\color{red} 0-2}
      }                     edge from parent node[left, draw=none] {\color{red} \ref*{def:ref-eq}-\ref*{rule:subs}}
      edge from parent node[right, draw=none] {\color{red} 0}
      }                     edge from parent node[left, draw=none] {\color{red} \ref*{rule:contradiction}}
      }
                          ;
\end{tikzpicture}
\end{center}

The next theorem allows to replace a reflexive equality by the unary attribute
of Definition \ref{def:ref-eq}.  This expresses a particular case of Rule \ref{rule:subs},
which is used in the proof of the theorem. In combination with Rule \ref{rule:contradiction}, the next
theorem allows to add the unary attribute for every object.

\begin{theorem}[Abbreviation reflexive equality, \ref{def:ref-eq}]\ \\ \label{thm:abb-ref-eq}
  $$\lall \xi=\xi\rall I\xi$$
\end{theorem}

  \begin{center}
\begin{tikzpicture}[every node/.style = {draw, shape=rectangle,anchor=north},growth parent anchor=south,level distance=.4 cm,sibling distance = 4.7cm]]
\node{$ \lexi \xi\eq \xi\rexi \ds{I}\xi $}
                    child{[level distance=.1cm]node{$x=x $}
                    child{[level distance=.4cm]node{$ \ds{I}x $}
                    child{node{$ Ix $}
          child{node{} 
edge from parent node[left, draw=none] {\color{red}  \ref*{rule:contradiction}}
edge from parent node[right, draw=none] {\color{red} 1}
      }
                     edge from parent node[left, draw=none] {\color{red} \ref*{def:ref-eq}-\ref*{rule:subs}}
      edge from parent node[right, draw=none] {\color{red} 1}
      }
                }
 edge from parent node[left, draw=none] {\color{red} \ref*{rule:exist}}
                    edge from parent node[right, draw=none] {\color{red} 0}
}
                    ;
\end{tikzpicture}
\end{center}

Every object is equal to itself as a consequence of Rule \ref{rule:contradiction}.
The following theorem expresses this fact. It states that if an object is activated with any
unary attribute, then it is equal to itself.

\begin{theorem}[Reflexive equality of arbitrary object]\ \\ \label{thm:reflexive-eq}
  $$\lall  \Gamma \xi\rall  \xi=\xi$$
\end{theorem}

  \begin{center}
\begin{tikzpicture}[every node/.style = {draw, shape=rectangle,anchor=north},growth parent anchor=south,level distance=.4 cm,sibling distance = 4.7cm]]
  \node{$ \lexi G \xi\rexi \xi\ds{\eq}\xi $}
                    child{[level distance=.1cm]node{$G x $}
                    child{[level distance=.1cm]node{$ x\ds{\eq} x $}
          child{node{} 
      }
      }
                     edge from parent node[left, draw=none] {\color{red} \ref*{rule:exist}}
                    edge from parent node[right, draw=none] {\color{red} 0}
                }
                    ;
\end{tikzpicture}
\end{center}

Given an object letter, that we assume activated with some unary attribute, we may assign it a new object
letter by means of an equality. This is expressed by the following theorem.

\begin{theorem}[New equal object, \ref{def:ref-eq}]\label{thm:new-name}\ \\
 $$\lall \Gamma \xi \rall \lexi I\eta \rexi \eta=\xi$$
\end{theorem}

\begin{center}
\begin{tikzpicture}[every node/.style = {draw, shape=rectangle,anchor=north},growth parent anchor=south,level distance=0.4cm, sibling distance = 6cm]
\node{$\lexi  G \xi \rexi\lall  I\eta \rall\eta \ds{\eq} \xi$} 
child{[level distance=0.1cm]
  node{$ G x$} 
          child{[level distance=0.4cm]node{$\lall I \eta \rall \eta \ds{\eq} x $} 
          child{node{$ Ix $} 
          child{[level distance=0.1cm]node{$x\ds{\eq} x$} 
          child{node{} 
      }
 edge from parent node[left, draw=none] {\color{red} \ref*{rule:all}}
edge from parent node[right, draw=none] {\color{red} 0-1}
      }
      edge from parent node[left, draw=none] {\color{red}
        \ref*{thm:abb-ref-eq}-\ref*{rule:contradiction}-\ref*{rule:exist}}
      }
    }
 edge from parent node[left, draw=none] {\color{red} \ref*{rule:exist}}
edge from parent node[right, draw=none] {\color{red} 0}
}
                    ;
\end{tikzpicture}
\end{center}

The attribute in the following definition
states that the object in question has an empty domain.
\begin {definition}[Attribute of empty domain]\label{def:empty}
$$O \xi \prop \lall \eta=\eta \rall \nonmem \eta \xi$$ 
\end{definition}

Theorem \ref{thm:zero-exists} shows existence of 
an object with empty domain. It is constructed by applying the rule of choice
with a condition that is impossible to satisfy.
 \begin{theorem}[Existence of empty domain, \ref{def:dom-prop}, \ref{def:empty}]\ \\ \label{thm:zero-exists}
$$\lexi O\xi \rexi {O}\xi$$ 
\end{theorem}
\begin{center}
		\begin{tikzpicture}[every node/.style = {draw, shape=rectangle,anchor=north},growth parent anchor=south,level distance=.1 cm,sibling distance = 5cm]
		\node{$\lall O\xi \rall \ds{O}\xi$} 
                child{[level distance = 0.1cm, sibling distance=3.5cm]
                  node{$\lall \ds{I}\xi \rall  \lexi I\eta \rexi \eta=\xi $}
          child{[level distance=.4cm]node{$\ds{M}(\ds{I})$}
          child{[level distance=.4cm]node{$y=y$}
          child{[level distance=.4cm]node{$\lexi \ds{I} \xi \rexi \nonmem \xi y$}
				child{node{$\ds{I} x$}
				child{[level distance = 0.1cm]node{$x\ds{\eq} x$}
				child{node{$ $}
                                }
         edge from parent node[left, draw=none] {\color{red} \ref*{def:ref-eq}-\ref*{rule:subs}}
         edge from parent node[right, draw=none] {\color{red} 0}
                                }
         edge from parent node[left, draw=none] {\color{red} \ref*{rule:exist}}
         edge from parent node[right, draw=none] {\color{red} 0}
                              }
         edge from parent node[left, draw=none] {\color{red} \ref*{def:dom-prop}-\ref*{rule:subs}}
         edge from parent node[right, draw=none] {\color{red} 0}
                                }
           edge from parent node[left, draw=none] {\color{red} \ref*{rule:contradiction}}
                                }
          }
				child{[level distance=.4cm]node{$M(\ds{I })y$}
				child{[sibling distance=3cm,level distance = 0.1cm]node{$D(\ds{I}) f$}
                                  child{[level distance = 0.4cm]node{$ {O}f$}
                                  child{[level distance = 0.4cm]node{$ \ds{O} f$}
                                  child{[level distance = 0.4cm]node{$ $}
         edge from parent node[right, draw=none] {\color{red} 1}
       }
         edge from parent node[left, draw=none] {\color{red} \ref*{rule:all}}
         edge from parent node[right, draw=none] {\color{red} 0-4}
       }
       }
                                    child{[level distance = 0.4cm]node{$ \ds{O}f$}
                                  child{[level distance = 0.4cm]node{$ \lexi \xi=\xi \rexi \mem \xi f$}
                                  child{[level distance = 0.1cm]node{$ x= x$}
                                  child{[level distance = 0.4cm]node{$ \mem x f$}
                                  child{node{$ \ds{I} x$}
                                  child{[level distance = 0.1cm]node{$ x\ds{\eq} x$}
                                  child{node{$ $}
       }
         edge from parent node[left, draw=none] {\color{red} \ref*{def:ref-eq}-\ref*{rule:subs}}
         edge from parent node[right, draw=none] {\color{red} 0}
       }
         edge from parent node[left, draw=none] {\color{red} \ref*{def:dom-prop}-\ref*{rule:subs}-\ref*{rule:exist}-\ref*{rule:all}}
         edge from parent node[right, draw=none] {\color{red} 0-4}
       }
     }
            edge from parent node[left, draw=none] {\color{red} \ref*{rule:exist}}
         edge from parent node[right, draw=none] {\color{red} 0}
}
          edge from parent node[left, draw=none] {\color{red} \ref*{def:empty}-\ref*{rule:subs}}
         edge from parent node[right, draw=none] {\color{red} 0}
}
       }
         edge from parent node[left, draw=none] {\color{red} \ref*{thm:choice}}
         edge from parent node[right, draw=none] {\color{red} 0-1}
       }}
     }
     child{[level distance = 0.4cm]node{$\lexi \ds{I} \xi\rexi  \lall I \eta\rall \eta\ds{\eq} \xi $}
       child{node{$\ds{I} x$}
       child{[level distance = 0.1cm]node{$x\ds{\eq} x$}
	child{node{$ $}
       }
         edge from parent node[left, draw=none] {\color{red} \ref*{def:ref-eq}-\ref*{rule:subs}}
         edge from parent node[right, draw=none] {\color{red} 0}
       }
         edge from parent node[left, draw=none] {\color{red} \ref*{rule:exist}}
         edge from parent node[right, draw=none] {\color{red} 0}
       }
       };
		\end{tikzpicture}
	\end{center}	

We use the above existence result to define an object with empty domain. We repeat the definition of the attribute of being empty, so that in future
references we do not have to refer to both Definitions \ref{def:empty} and \ref{def:empty-object} in the same preamble.

\begin{definition} [Object with empty domain] \ \\ \label{def:empty-object}
\end{definition}
\begin{center}
\begin{tikzpicture}[every node/.style = {draw, shape=rectangle,anchor=north},growth parent anchor=south,level distance=.4 cm,sibling distance = 4.7cm]]
  \node{$O\xi \prop O_{\ref{def:empty}}\xi$}
  child{[level distance=.1cm]node{$Oo$}
        edge from parent node[left, draw=none] {\color{red} \ref*{thm:zero-exists}}
  };
    \end{tikzpicture}
\end{center}

The following theorem expresses uniqueness of the object with an empty domain. The proof is
an application of Theorem \ref{thm:function-uniqueness}.

 \begin{theorem}[Uniqueness of empty domain, \ref{def:empty-object}]\ \\ \label{thm:zero-unique}
$$\lall O \xi \rall \xi= o$$ 
\end{theorem}
\begin{center}
		\begin{tikzpicture}[every node/.style = {draw, shape=rectangle,anchor=north},growth parent anchor=south,level distance=.4 cm,sibling distance = 6.4cm]
                  \node{$\lexi O \xi \rexi \xi\ds{\eq} o$} 
                child{[level distance=.1cm]
                  node{$O f$}
                child{[level distance=.4cm]
                  node{$f \ds{\eq} o$}
                child{[level distance=.4cm]
                  node{$\lall \xi=\xi \rall \nonmem \xi f$}
     child{[level distance = 0.4cm]node{$\lexi f (\xi)\ds{\eq} o(\xi)\rexi  \lall \nonmem
                    \xi f \rall \mem \xi o $}
       child{[level distance=.1cm]node{$ f (x)\ds{\eq} o (x)$}
       child{[level distance=.4cm]node{$\lall \nonmem x f \rall \mem x o$}
	child{node{$\nonmem x f $}
	child{node{$\mem x o$}
	child{node{$O o$}
	child{node{$\nonmem x o$}
	child{node{$$}
         edge from parent node[right, draw=none] {\color{red} 2}
        }
        edge from parent node[left, draw=none] {\color{red} \ref*{def:empty-object}-\ref*{def:empty}-
          \ref*{rule:contradiction}-\ref*{rule:all}}
         edge from parent node[right, draw=none] {\color{red} 0}
       }
        edge from parent node[left, draw=none] {\color{red} \ref*{def:empty-object}}
       }
          edge from parent node[left, draw=none] {\color{red} \ref*{rule:all}}
         edge from parent node[right, draw=none] {\color{red} 0-1}
     }
          edge from parent node[left, draw=none] {\color{red} \ref*{rule:all}-\ref*{rule:contradiction}}
         edge from parent node[right, draw=none] {\color{red} 3}
  }
        }
          edge from parent node[left, draw=none] {\color{red} \ref*{rule:exist}}
         edge from parent node[right, draw=none] {\color{red} 0}
       }
     edge from parent node[left, draw=none] {\color{red} \ref*{thm:function-uniqueness}-\ref*{rule:all}-\ref*{rule:contradiction}
     }
         edge from parent node[right, draw=none] {\color{red} 1}
     }
   edge from parent node[left, draw=none] {\color{red} \ref*{def:empty-object}-\ref*{def:empty}
     \ref*{rule:subs}}
         edge from parent node[right, draw=none] {\color{red} 1}
       }}
         edge from parent node[left, draw=none] {\color{red} \ref*{rule:exist}}
         edge from parent node[right, draw=none] {\color{red} 0}
       };
		\end{tikzpicture}
	\end{center}	

Our next goal is to find an object with nonempty domain. To this end, the following theorem
expresses existence  of an object whose domain contains the object with empty domain.
We use Theorem \ref{thm:universe}.

 \begin{theorem}[Non-empty domain, \ref{def:ref-eq}, \ref{def:empty-object}]\ \\ \label{thm:non-empty-domain}
$$\lexi I\xi \rexi \mem o \xi $$ 
\end{theorem}
\begin{center}
		\begin{tikzpicture}[every node/.style = {draw, shape=rectangle,anchor=north},growth parent anchor=south,level distance=.4 cm,sibling distance = 6.4cm]
                  \node{$\lall I\xi \rall \nonmem o \xi$} 
                child{[level distance=.4cm]
                  node{$\mem o f$}
                child{[level distance=.4cm]
                  node{$If$}
               child{[level distance=.4cm]
                  node{$\nonmem o f$}
	child{node{$ $}
         edge from parent node[right, draw=none] {\color{red} 2}
  }
     edge from parent node[left, draw=none] {\color{red} \ref*{rule:all}}
         edge from parent node[right, draw=none] {\color{red} 0-2}
       }
     edge from parent node[left, draw=none] {\color{red} \ref*{rule:contradiction}-\ref*{thm:abb-ref-eq}}
       }
   edge from parent node[left, draw=none] {\color{red} 
       \ref*{thm:universe}-\ref*{rule:contradiction}-\ref*{rule:all}-\ref*{rule:exist}}
       };
		\end{tikzpicture}
	\end{center}	

        Having both an object with empty domain and an object with non-empty domain, we are ready
        to show that for any object one can find an object that is different.
        This theorem is somewhat parallel to Theorem \ref{thm:new-name}.

 \begin{theorem}[Different object, \ref{def:ref-eq}, \ref{def:empty-object}]\ \\ \label{thm:different-object}
$$\lall \Gamma \xi\rall \lexi I \eta\rexi \eta\ds{\eq} \xi $$ 
\end{theorem}
\begin{center}
\begin{tikzpicture}[every node/.style = {draw, shape=rectangle,anchor=north},growth parent anchor=south,level distance=0.4cm, sibling distance = 6cm]
\node{$\lexi  G \xi \rexi\lall  I \eta\rall\eta {\eq} \xi$} 
child{[level distance=0.1cm]
  node{$ G x$} 
          child{[level distance=0.1cm]node{$\lall I \eta\rall \eta {\eq} x$} 
          child{[level distance=0.4cm]node{$o \ds{\eq} x$} 
          child{node{$Io$} 
          child{node{$o {\eq} x$} 
          child{node{} 
edge from parent node[right, draw=none] {\color{red} 2}
      }
      edge from parent node[left, draw=none] {\color{red} \ref*{rule:all}}
      edge from parent node[right, draw=none] {\color{red} 0-2}
      }
      edge from parent node[left, draw=none] {\color{red} \ref*{rule:contradiction}-\ref*{thm:abb-ref-eq}}
      }
      }
          child{[level distance=0.4cm]node{$o \eq x$} 
          child{[level distance=0.1cm]node{$Iy$} 
            child{[level distance=0.4cm]node{$\mem o y$}
                        child{node{$ y\eq x$} 
          child{node{$y\eq o$} 
          child{node{$\mem o o$} 
          child{node{$\nonmem o o$} 
          child{node{$ $} 
edge from parent node[right, draw=none] {\color{red} 1}
      }
      edge from parent node[left, draw=none] {\color{red} \ref*{def:empty-object}-\ref*{def:empty}-\ref*{rule:subs}-\ref*{rule:all}-\ref*{rule:contradiction}}
      }
 edge from parent node[left, draw=none] {\color{red} \ref*{rule:subf}}
edge from parent node[right, draw=none] {\color{red} 0-2}
      }
 edge from parent node[left, draw=none] {\color{red} \ref*{rule:subf}}
edge from parent node[right, draw=none] {\color{red} 0-3}
      }
 edge from parent node[left, draw=none] {\color{red} \ref*{rule:all}}
edge from parent node[right, draw=none] {\color{red} 1-3}
      }
      }
 edge from parent node[left, draw=none] {\color{red} \ref*{thm:non-empty-domain}-\ref*{rule:exist}}
      }
      }
  }
 edge from parent node[left, draw=none] {\color{red} \ref*{rule:exist}}
edge from parent node[right, draw=none] {\color{red} 0}
}
                    ;
\end{tikzpicture}
\end{center}

The following two theorems together constitute Russell's paradox.
The first theorem states that for every object there
is another object  that is either equal or not equal to the first object depending on
a certain condition on the first object. 
The proof prominently uses Theorem \ref{thm:different-object} to produce the different object
when needed.

The second theorem states that this second object  cannot
be chosen prior to the first object, even if we allow
it to float with the first object. 
If we have choose the second object before the first,
we may then choose the first object be equal to the second object to obtain a contradiction.

\begin{theorem}[Russell's paradox, first part, \ref{def:ref-eq}]\ \\
\label{thm:russell}
$$
\lall  I \xi \rall \lexi  I \eta \rexi 
\lexi  \lall  \xi(\xi)\eq \xi\rall \eta\ds{\eq} \xi \rexi 
\lall  \xi(\xi)\ds{\eq} \xi\rall \eta {\eq}  \xi
$$

\end{theorem}

  \begin{center}
\begin{tikzpicture}[every node/.style = {draw, shape=rectangle,anchor=north},growth parent anchor=south,level distance=0.4cm, sibling distance = 4cm]
\node{$\lexi I \xi \rexi \lall  I \eta \rall 
\lall  \lall  \xi(\xi)\eq \xi\rall \eta\ds{\eq} \xi \rall 
\lexi  \xi(\xi)\ds{\eq} \xi\rexi \eta\ds{\eq}  \xi$} 
child{[level distance=0.1cm]node{$ I x$} 
          child{[sibling distance = 5.5cm]node{$ \lall  I \eta \rall 
\lall  \lall  x(x)\eq x\rall \eta\ds{\eq} x \rall 
\lexi  x(x)\ds{\eq} x\rexi \eta\ds{\eq} x $ } 
          child{[level distance=0.4cm]node{$x(x)\eq x$} 
            child{[level distance=0.1cm]node{$I y$} 
            child{[level distance=0.4cm]node{$y\ds{\eq} x$} 
      child{[level distance=0.1cm, sibling distance = 3.3cm]node
{$\lall  \lall  x(x)\eq x\rall y\ds{\eq} x \rall \lexi  x(x)\ds{\eq} x\rexi y\ds{\eq} x $} 
            child{[level distance=0.4cm]node{$\lall x(x)\eq  x\rall y\ds{\eq} x$}
            child{node{$x(x)\ds{\eq} x$}
          child{node{} 
edge from parent node[right, draw=none] {\color{red} 5}
}
 edge from parent node[left, draw=none] {\color{red} \ref*{rule:all}-\ref*{rule:exist}}
edge from parent node[right, draw=none] {\color{red} 0-1}
}
}
             child{[level distance=0.4cm]node{$\lexi x(x)\eq  x\rexi y\eq  x$} 
              child{[level distance=0.4cm]node{$y\eq  x$} 
          child{node{} 
edge from parent node[right, draw=none] {\color{red} 3}
}
        edge from parent node[left, draw=none] {\color{red} \ref*{rule:exist}}
        edge from parent node[right, draw=none] {\color{red} 0}
}
      }
        edge from parent node[left, draw=none] {\color{red} \ref*{rule:exist}}
        edge from parent node[right, draw=none] {\color{red} 1-3}
}
  }
         edge from parent node[left, draw=none] {\color{red} \ref*{thm:different-object}-\ref*{rule:all}-\ref*{rule:exist}}
        edge from parent node[right, draw=none] {\color{red} 2}
           }
         }
          child{[level distance=0.4cm]node{$x(x)\ds{\eq} x$} 
            child{[level distance=0.1cm,sibling distance = 3.3cm]node
{$\lall  \lall  x(x)\eq x\rall x\ds{\eq} x \rall \lexi  x(x)\ds{\eq} x\rexi x\ds{\eq} x $} 
          child{[level distance=2.1cm]node{$\lexi x(x)\eq x\rexi x\eq x$}
          child{[level distance=0.4cm]node{$ x(x)\eq x$}
          child{node{} 
edge from parent node[right, draw=none] {\color{red} 3}
}
edge from parent node[left, draw=none] {\color{red} \ref*{rule:exist}}
 edge from parent node[right, draw=none] {\color{red} 0}    
}}
          child{[level distance=0.4cm]node{$\lall x(x)\eq x \rall x\ds{\eq} x$}
          child{[level distance=0.1cm]node{$x\ds{\eq} x$}
          child{node{} 
}
       edge from parent node[left, draw=none] {\color{red} \ref*{rule:all}-\ref*{rule:exist}}
      edge from parent node[right, draw=none] {\color{red} 0-1}
  }
      } 
         edge from parent node[left, draw=none] {\color{red} \ref*{rule:all}}
        edge from parent node[right, draw=none] {\color{red} 1-2}
    }}}
         edge from parent node[left, draw=none] {\color{red} \ref*{rule:exist}}
        edge from parent node[right, draw=none] {\color{red} 0}
    }
                    ;
\end{tikzpicture}
\end{center}


\begin{theorem}[Russell's paradox, second part, \ref{def:ref-eq}]\ \\
\label{thm:russell2}
$$\lall  I \eta \rall \lexi  I \xi \rexi 
\lall  \lall  \xi(\xi)\eq \xi\rall \eta(\xi)\ds{\eq} \xi \rall 
\lexi  \xi(\xi)\ds{\eq} \xi\rexi \eta(\xi)\ds{\eq}   \xi$$
\end{theorem}

\begin{center}
\begin{tikzpicture}[every node/.style = {draw, shape=rectangle,anchor=north},growth parent anchor=south,level distance=0.4cm, sibling distance = 4cm]
\node{$\lexi  I \eta \rexi \lall  I \xi \rall 
\lexi  \lall  \xi(\xi)\eq \xi\rall \eta(\xi)\ds{\eq} \xi \rexi 
\lall  \xi(\xi)\ds{\eq} \xi\rall \eta(\xi)\eq   \xi$} 
child{[level distance=0.1cm]node{$ I y$} 
  child{[level distance=0.4cm]node{$ \lall  I \xi \rall 
\lexi  \lall  \xi(\xi)\eq \xi\rall y(\xi)\ds{\eq} \xi \rexi 
\lall  \xi(\xi)\ds{\eq} \xi\rall y(\xi)\eq  \xi $ 
} 
          child{[sibling distance = 7cm]node{
$  
\lexi  \lall  y(y)\eq y\rall y(y)\ds{\eq} y \rexi
\lall  y(y)\ds{\eq} y\rall y(y)\eq  y $ 
} 
          child{[level distance=0.1cm]node{$  \lall  y(y)\eq y\rall y(y)\ds{\eq} y $} 
            child{node{$\lall y(y)\ds{\eq} y\rall y(y)\eq y $} 
            child{[level distance=0.4cm]node{$y(y)\ds{\eq} y$} 
           child{node
             {$y(y)\eq y$}
                       child{node{} 
edge from parent node[right, draw=none] {\color{red} 1}
      }
                        edge from parent node[left, draw=none] {\color{red} \ref*{rule:all}}
                       edge from parent node[right, draw=none] {\color{red} 0-1}
                    }}
          child{[level distance=0.4cm]node{$ y(y)\eq y$}
          child{node{$ y(y)\ds{\eq} y$}
          child{node{} 
edge from parent node[right, draw=none] {\color{red} 1}
}
edge from parent node[left, draw=none] {\color{red} \ref*{rule:all}}
             edge from parent node[right, draw=none] {\color{red} 0-2}
         }}}
          edge from parent node[left, draw=none] {\color{red} \ref*{rule:exist}}
         edge from parent node[right, draw=none] {\color{red} 0}
     }
          edge from parent node[left, draw=none] {\color{red} \ref*{rule:all}}
         edge from parent node[right, draw=none] {\color{red} 0-1}
     }}
          edge from parent node[left, draw=none] {\color{red} \ref*{rule:exist}}
         edge from parent node[right, draw=none] {\color{red} 0}
     }
                    ;
\end{tikzpicture}
\end{center}

The term in $M(\Gamma)$ in Definition \ref{def:dom-prop} states that the attribute $\Gamma$
is controlled by the domain of a function. The next theorem shows that this sub domain property is
inherited by any more restrictive attribute.

\begin{theorem}[Nested attributes, \ref{def:dom-prop}]\ \\
  \label{thm:nested-att}
$$ \lall \lall \Gamma \xi \rall \Delta \xi \rall \lall M(\Delta)\rall M(\Gamma) $$ 
\end{theorem}
\begin{center}
		\begin{tikzpicture}[every node/.style = {draw, shape=rectangle,anchor=north},growth parent anchor=south,level distance=.4 cm,sibling distance = 6.4cm]
                  \node{$ \lexi \lall G \xi \rall D \xi \rexi \lexi M(D)\rexi \ds M(G) $} 
                child{[level distance=.1cm]node{$\lall G\xi \rall D \xi$}
                child{[level distance=.1cm]node{$M(D) $}
                child{[level distance=.4cm]node{$\ds M (G)$}
       child{[level distance=.1cm]node{$f=f$}
	child{[level distance=.4cm]node{$ \lall D \xi \rall \mem \xi f$}
	child{[level distance=.1cm]node{$ G x$}
	child{[level distance=.4cm]node{$ \nonmem x f$}
	child{node{$ D x$}
	child{node{$ \mem x f$}
	child{node{$ $}
         edge from parent node[right, draw=none] {\color{red} 2}
       }
               edge from parent node[left, draw=none] {\color{red} \ref*{rule:all}}
         edge from parent node[right, draw=none] {\color{red} 0-3}
       }
               edge from parent node[left, draw=none] {\color{red} \ref*{rule:all}}
         edge from parent node[right, draw=none] {\color{red} 1-6}
       }
       }
               edge from parent node[left, draw=none] {\color{red} \ref*{def:dom-prop}-\ref*{rule:subs}-\ref*{rule:all}-\ref*{rule:exist}}
         edge from parent node[right, draw=none] {\color{red} 1-2}
       }
     }           
         edge from parent node[left, draw=none] {\color{red} \ref*{def:dom-prop}-\ref*{rule:subs}-\ref*{rule:exist}}
         edge from parent node[right, draw=none] {\color{red} 1}
       }
       }
     }
         edge from parent node[left, draw=none] {\color{red} \ref*{rule:exist}}
         edge from parent node[right, draw=none] {\color{red} 0}
       };
		\end{tikzpicture}
	\end{center}	

The abbreviation $C\phi$ in Definition \ref{def:id-ext-succ}
states that $\phi$ is an identity object, meaning for each object
in its domain, it  floats to become this very object.

\begin{theorem}[Empty object as identity object, \ref{def:id-ext-succ}, \ref{def:empty-object}]\ \\
  \label{thm:empty-id}
$$ Co$$ 
\end{theorem}
\begin{center}
		\begin{tikzpicture}[every node/.style = {draw, shape=rectangle,anchor=north},growth parent anchor=south,level distance=.4 cm,sibling distance = 6.4cm]
                  \node{$\ds {C} o$} 
                child{[level distance=.4cm]node{$\lexi \mem \xi o \rexi o(\xi)\ds{\eq} \xi $}
                child{[level distance=.4cm]node{$\mem x o $}
                child{[level distance=.4cm]node{$\lall \xi\eq  \xi  \rall \nonmem \xi o$}
       child{[level distance=.4cm]node{$ \nonmem x o$}
	child{node{$ $}
         edge from parent node[right, draw=none] {\color{red} 2}
       }
     edge from parent node[left, draw=none] {\color{red} \ref*{rule:contradiction}-\ref*{rule:all}}
       edge from parent node[right, draw=none] {\color{red} 0}
     }
   edge from parent node[left, draw=none] {\color{red} \ref*{def:empty}-\ref*{def:empty-object}-\ref*{rule:subs}}
       }
     edge from parent node[left, draw=none] {\color{red} \ref*{rule:exist}}
        edge from parent node[right, draw=none] {\color{red} 0}
     }
         edge from parent node[left, draw=none] {\color{red} \ref*{def:id-ext-succ}-\ref*{rule:subs}}
         edge from parent node[right, draw=none] {\color{red} 0}
       };
		\end{tikzpicture}
	\end{center}	

The following theorem states that for every sub domain attribute, there is an identity object with  domain precisely described by this attribute.

        \begin{theorem}[Existence identity function, \ref{def:dom-prop}, \ref{def:id-ext-succ}]\
          \\ \label{thm:exist-identity}
$$\lall M(\Gamma) \rall \lexi D(\Gamma)\psi \rexi C\psi$$ 
\end{theorem}
\begin{center}
		\begin{tikzpicture}[every node/.style = {draw, shape=rectangle,anchor=north},growth parent anchor=south,level distance=.4 cm,sibling distance = 6.4cm]
                  \node{$\lexi M(G) \rexi \lall D(G)\psi \rall \ds{C}\psi$} 
                child{[level distance=.1cm]
                  node{$M(G)$}
                child{[level distance=.4cm]
                  node{$\lall D(G)\psi \rall \ds{C}\psi$}
                child{[level distance=.4cm]
                  node{$\lall G\xi \rall \lexi I\eta\rexi \eta=\xi$}
       child{[level distance=.1cm]node{$ D(G) g$}
       child{[level distance=.4cm]node{$ \lall G\xi \rall \lexi Ig(\xi)\rexi g(\xi)=\xi$}
	child{node{$\ds{C} g $}
	child{node{$\lexi \mem  \xi g \rexi g(\xi)\ds{\eq} \xi$}
	child{[level distance=.1cm]node{$\mem x g$}
	child{[level distance=.4cm]node{$g(x)\ds{\eq} x$}
	child{node{$\lall \mem \xi g \rall G\xi$}
	child{node{$Gx$}
	child{node{$g(x)=x$}
	child{node{$ $}
         edge from parent node[right, draw=none] {\color{red} 3}
        }
         edge from parent node[left, draw=none] {\color{red} \ref*{rule:all}-\ref*{rule:exist}}
         edge from parent node[right, draw=none] {\color{red} 0-6}
        }
          edge from parent node[left, draw=none] {\color{red} \ref*{rule:all}}
         edge from parent node[right, draw=none] {\color{red} 0-2}
        }
          edge from parent node[left, draw=none] {\color{red} \ref*{def:dom-prop}-\ref*{rule:subs}-\ref*{rule:exist}}
         edge from parent node[right, draw=none] {\color{red} 5}
        }
       }
          edge from parent node[left, draw=none] {\color{red} \ref*{rule:exist}}
         edge from parent node[right, draw=none] {\color{red} 0}
     }
          edge from parent node[left, draw=none] {\color{red} \ref*{def:id-ext-succ}-\ref*{rule:subs}}
         edge from parent node[right, draw=none] {\color{red} 0}
  }
         edge from parent node[left, draw=none] {\color{red} \ref*{rule:all}}
         edge from parent node[right, draw=none] {\color{red} 1-3}
        }
       }
     edge from parent node[left, draw=none] {\color{red} \ref*{thm:choice}}
       edge from parent node[right, draw=none] {\color{red} 0-2}
     }
   edge from parent node[left, draw=none] {\color{red} \ref*{thm:new-name}}
       }
     }
         edge from parent node[left, draw=none] {\color{red} \ref*{rule:exist}}
         edge from parent node[right, draw=none] {\color{red} 0}
       };
		\end{tikzpicture}
	\end{center}	

A successor of an object is an identity object whose domain contains precisely the original object as well as all objects in
the domain of the original object. The attribute $E\phi\xi$ states that $\xi$ is in this extended domain
of $\phi$. The attribute $H\psi\phi$ states that $\psi$ is the successor of $\phi$. 

The purpose of the next two theorems is to show without using Theorem \ref{thm:infinity} that every identity object has a successor.
 The existence of such a successor is part of Theorem \ref{thm:infinity}, but the main point of Theorem \ref{thm:infinity}
 is the existence of some object whose domain contains a successor of each of the objects in its domain.
 
We first show with the aid of Theorem \ref{thm:universe} that for every identity object there is an object whose domain is at least the extended object of the original domain. 

\begin{theorem}[Extension of domain, \ref{def:range-union-power}, \ref{def:id-ext-succ}]\ \\  \label{thm:domainextension}
$$ \lall C\phi\rall \lexi\psi=\psi\rexi \lall E\phi \xi\rall \mem \xi \psi $$
\end{theorem}  
	\begin{center}
		\begin{tikzpicture}[every node/.style = {draw, shape=rectangle,anchor=north},growth parent anchor=south,level distance=.4 cm,sibling distance = 5.3cm]
                  \node{$ \lexi C\phi\rexi \lall\psi=\psi\rall\lexi E\phi \xi\rexi \nonmem \xi \psi
                    $}
                  	child{[level distance=.1cm]node{$Cf$} 
		child{[level distance=.4cm]node{$
				\lall \psi=\psi \rall 
				\lexi Ef \xi\rexi \nonmem \xi \psi
				$} 
			child{[level distance=.1cm]node{$\mem f g$} 
		child{[level distance=.4cm]node{$R f g$}
                                          child{[level distance=.4cm]node{$\lall \mem \xi f \rall \mem {f(\xi)} g$}
					 			child{[level distance=.1cm]node{$ E f x$}
											child{[level distance=.4cm]node{$ \nonmem x g $}						
												child{[level distance=.1cm]node{$ \lall \mem \xi f \rall f(\xi)=\xi$}	
												child{[level distance=.4cm]node{$ \nonmem x f$}	
												child{node{$ \lall \nonmem x f \rall f=x$}	
												child{node{$ f=x$}	
												child{node{$ \nonmem f g $}	
												child{node{$ $}	
         edge from parent node[right, draw=none] {\color{red} 9}
                                                                                                }
         edge from parent node[left, draw=none] {\color{red} \ref*{rule:subf}}
         edge from parent node[right, draw=none] {\color{red} 0-4}
                                                                                                }
       edge from parent node[left, draw=none] {\color{red} \ref*{rule:all}}
         edge from parent node[right, draw=none] {\color{red} 0-1}
                                                                                                }
         edge from parent node[left, draw=none] {\color{red} \ref*{def:id-ext-succ}-\ref*{rule:subs}}
         edge from parent node[right, draw=none] {\color{red} 3}
                                                                                                }
       }
        child{[level distance=.4cm]node{$ \mem x f$}
												child{node{$ f(x)=x$}	
												child{node{$ \mem {f(x)} g$}	
												child{node{$ \mem x g$}	
												child{node{$ $}	
         edge from parent node[right, draw=none] {\color{red} 5}
                                                                                                }
           edge from parent node[left, draw=none] {\color{red} \ref*{rule:subf}}
         edge from parent node[right, draw=none] {\color{red} 0-1}
                                                                                                }
           edge from parent node[left, draw=none] {\color{red} \ref*{rule:all}}
         edge from parent node[right, draw=none] {\color{red} 1-5}
                                                                                                }
       edge from parent node[left, draw=none] {\color{red} \ref*{rule:all}}
         edge from parent node[right, draw=none] {\color{red} 0-1}
                                                                                                }
                                                                                                }
         edge from parent node[left, draw=none] {\color{red} \ref*{def:id-ext-succ}-\ref*{rule:subs}}
         edge from parent node[right, draw=none] {\color{red} 6}                                                                                      }
       }
    edge from parent node[left, draw=none] {\color{red} \ref*{rule:contradiction}-\ref*{rule:all}-\ref*{rule:exist}}
         edge from parent node[right, draw=none] {\color{red} 3}
       }
       edge from parent node[left, draw=none] {\color{red} \ref*{def:range-union-power}-\ref*{rule:subs}}
         edge from parent node[right, draw=none] {\color{red} 0}
     }
       }
              edge from parent node[left, draw=none] {\color{red} \ref*{thm:universe}-\ref*{rule:contradiction}-\ref*{rule:exist}}
       }
}
       edge from parent node[left, draw=none] {\color{red} \ref*{rule:exist}}
         edge from parent node[right, draw=none] {\color{red} 0}
};
		\end{tikzpicture}
	\end{center}	

The next theorem uses the previous two theorems  to show that every identity object has a successor.   

\begin{theorem}[Existence of successor, \ref{def:dom-prop}, \ref{def:id-ext-succ}]\ \\ \label{thm:successor-existence}
$$\lall C\phi\rall \lexi \psi=\psi \rexi {H}\psi \phi$$
\end{theorem}

\begin{center}
		\begin{tikzpicture}[every node/.style = {draw, shape=rectangle,anchor=north},growth parent anchor=south,level distance=.4 cm,sibling distance = 5.4cm]
                  \node{$\lexi C\phi \rexi \lall \psi=\psi \rall
                        \ds{H}\psi \phi$} 
                child{[level distance=.1cm]
                  node{$Cf$}
                child{[level distance=.4cm]
                  node{$\lall \psi=\psi \rall
                        \ds{H}\psi f $}
                child{[level distance=.1cm]
                  node{$h=h$}
     child{[level distance = 0.4cm]node{$\lall Ef\xi\rall \mem \xi h $}
       child{[level distance=.1cm]node{$ A\xi:Ef\xi$}
	child{[level distance=.4cm]node{$\ds{M}(A)$}
	child{[level distance=.4cm]node{$  \lall \phi=\phi  \rall \lexi A\xi\rexi \nonmem \xi \phi$}
	child{[level distance=.4cm]node{$ \lexi A\xi \rexi \nonmem \xi h$}
	child{node{$ \lexi Ef\xi \rexi \nonmem \xi h$}
	child{node{$ $}
         edge from parent node[right, draw=none] {\color{red} 5}
        }
          edge from parent node[left, draw=none] {\color{red} \ref*{rule:subs}}
         edge from parent node[right, draw=none] {\color{red} 0-3}
        }
          edge from parent node[left, draw=none] {\color{red} \ref*{rule:all}}
         edge from parent node[right, draw=none] {\color{red} 0-4}
        }
          edge from parent node[left, draw=none] {\color{red} \ref*{def:dom-prop}-\ref*{rule:subs}}
         edge from parent node[right, draw=none] {\color{red} 0}
        }
        }	child{[level distance=.4cm]node{$ M(A)$}
	child{[level distance=.1cm]node{$ D(A)g$}
	child{[level distance=.4cm]node{$ Cg$}
	child{node{$ \ds{H} gf$}
          child{node{$ \lall \lall Ef\xi \rall \mem \xi g\rall
              \lall \lall \mem \xi g  \rall Ef\xi\rall
              \ds{C}g
              $}
	child{node{$ \lexi \lall A\xi \rall \mem \xi g\rexi
              \lall \mem \xi g  \rall A\xi$}
	child{node{$ \lexi \lall Ef\xi \rall \mem \xi g\rexi
              \lall \mem \xi g  \rall Ef\xi$}
	child{[level distance=.1cm]node{$  \lall Ef\xi \rall \mem \xi g$}
	child{[level distance=.4cm]node{$ \lall \mem \xi g  \rall Ef\xi$}
	child{node{$ \ds{C}g$}
	child{node{$ $}
         edge from parent node[right, draw=none] {\color{red} 7}
        }
          edge from parent node[left, draw=none] {\color{red} \ref*{rule:all}}
         edge from parent node[right, draw=none] {\color{red} 0-1-4}
        }
        }
          edge from parent node[left, draw=none] {\color{red} \ref*{rule:exist}}
         edge from parent node[right, draw=none] {\color{red} 0}
        }
          edge from parent node[left, draw=none] {\color{red} \ref*{rule:subs}}
         edge from parent node[right, draw=none] {\color{red} 0-6}
        }
          edge from parent node[left, draw=none] {\color{red} \ref*{def:dom-prop}-\ref*{rule:subs}}
         edge from parent node[right, draw=none] {\color{red} 3}
        }
          edge from parent node[left, draw=none] {\color{red} \ref*{def:id-ext-succ}-\ref*{rule:subs}}
         edge from parent node[right, draw=none] {\color{red} 0}
        }
         edge from parent node[left, draw=none] {\color{red} \ref*{rule:contradiction}-\ref*{rule:all}}
         edge from parent node[right, draw=none] {\color{red} 6}
        }
       }
          edge from parent node[left, draw=none] {\color{red} \ref*{thm:exist-identity}-\ref*{rule:exist}}
         edge from parent node[right, draw=none] {\color{red} 0}
     }
  }
      edge from parent node[left, draw=none] {\color{red} \ref*{rule:con-pro}}
       }
     }
   edge from parent node[left, draw=none] {\color{red} \ref*{thm:domainextension}-\ref*{rule:exist}}
        edge from parent node[right, draw=none] {\color{red} 1}
       }}
         edge from parent node[left, draw=none] {\color{red} \ref*{rule:exist}}
         edge from parent node[right, draw=none] {\color{red} 0}
       };
		\end{tikzpicture}
	\end{center}	

The next theorem states that an identity object can not be in its own domain.
In particular, for an identity object the extended domain is strictly larger than the domain of the object.
Beginning with the empty object, which is an identity object, one may therefore use the previous
theorem  to construct larger and larger  finite domains.  
As discussed more thoroughly in the next section, by Theorem \ref{thm:infinity}, all these objects can be assumed
in  the domain of some object, which therefore has an infinite domain.

\begin{theorem}[Identity function not in its domain, \ref{def:id-ext-succ}]\ \\
  \label{thm:non-self-id}
$$\lall C\phi \rall \nonmem \phi \phi$$ 
\end{theorem}
\begin{center}
		\begin{tikzpicture}[every node/.style = {draw, shape=rectangle,anchor=north},growth parent anchor=south,level distance=.4 cm,sibling distance = 6.4cm]
                  \node{$\lexi C\phi \rexi \mem \phi \phi$} 
                child{[level distance=.1cm]
                  node{$Cf $}
                child{[level distance=.4cm]
                  node{$\mem f f $}
                child{[level distance=.4cm]
                  node{$\lall \mem \xi f \rall f(\xi) {\eq} \xi$}
       child{[level distance=.4cm]node{$ f(f)\eq f$}
	child{node{$ $}
         edge from parent node[right, draw=none] {\color{red} 2}
       }
     edge from parent node[left, draw=none] {\color{red} \ref*{rule:all}}
       edge from parent node[right, draw=none] {\color{red} 0-1}
     }
   edge from parent node[left, draw=none] {\color{red} \ref*{def:id-ext-succ}-\ref*{rule:subs}}
         edge from parent node[right, draw=none] {\color{red} 1}
       }
     }
         edge from parent node[left, draw=none] {\color{red} \ref*{rule:exist}}
         edge from parent node[right, draw=none] {\color{red} 0}
       };
		\end{tikzpicture}
	\end{center}	

\section{Natural numbers}
\label{section:peano}

In this chapter we discuss the Peano axioms of the natural
numbers. We will construct a model of the natural numbers and then
prove that this model satisfies the Peano axioms.

We use the object with empty domain as zero of the natural numbers.
An inductive domain contains zero and with each of its objects also
a successor of the object, with successor being the identity object on the extended domain
as described in Definition \ref{def:id-ext-succ}. The attribute $S$ in the following definition
states that $\phi$ has an inductive domain. An object that is in every inductive domain
is called a natural number, this is expressed by the attribute $N$.
The attribute $T$ describes an object with inductive domain consisting entirely  of natural numbers.

\begin{definition}[Inductive domains, natural numbers, \ref{def:dom-prop}, \ref{def:id-ext-succ}, \ref{def:empty-object}]\label{def:ind-nat}
\end{definition}
\begin{center}
\begin{tikzpicture}[every node/.style = {draw, shape=rectangle,anchor=north},growth parent anchor=south,level distance=.4 cm,sibling distance = 4.7cm]]
  \node{$S\phi \prop \lexi \mem o \phi \rexi \lall \mem \xi \phi \rall \lexi \mem \eta \phi \rexi H\eta \xi$}
  child{[level distance=.1cm]node{$N \xi \prop \lall S\phi \rall  \mem \xi \phi$}
  child{[level distance=.1cm]node{$T\sigma \prop  \lexi D(N) \sigma  \rexi \lall N\xi \rall H\sigma (\xi)\xi$}
}
  };
    \end{tikzpicture}
\end{center}

The following theorem shows existence of an inductive domain. The proof uses the rule of infinity.
  \begin{theorem}[Existence of inductive object,  \ref{def:id-ext-succ}, \ref{def:empty-object}, \ref{def:ind-nat}]\ \\ \label{thm:succ-funct}
    $$\lexi \phi =\phi \rexi  {S}\phi $$
  \end{theorem}
  
  \begin{center}
		\begin{tikzpicture}[every node/.style = {draw, shape=rectangle,anchor=north},growth parent anchor=south,level distance=.4 cm,sibling distance = 4.7cm]
		\node{$\lall \phi =\phi \rall  \ds{S}\phi $}
                                  child{[level distance=.4cm]node{$\lexi \mem o s \rexi \lall \mem \xi s \rall \lexi \mem \eta s \rexi H\eta \xi$}
                                  child{[level distance=.4cm]node{$\ds{S} s $}
				child{[level distance=.4cm]node{$\lall \mem 0 s\rall \lexi \mem \xi s \rexi \lall \mem \eta s\rall  \ds{H} \eta \xi $}
                                  child{[level distance=.4cm]node{$ $}
edge from parent node[right, draw=none] {\color{red} 2}
                                }
                                    edge from parent node[left, draw=none] {\color{red} \ref*{def:ind-nat}-\ref*{rule:subs}}
edge from parent node[right, draw=none] {\color{red} 0}
                          }
                                    edge from parent node[left, draw=none] {\color{red} \ref*{rule:contradiction}-\ref*{rule:all}}
edge from parent node[right, draw=none] {\color{red} 1}
                        }
                        edge from parent node[left, draw=none] {\color{red} \ref*{thm:infinity}-\ref*{rule:contradiction}-\ref*{rule:all}}
edge from parent node[right, draw=none] {\color{red} 0}
                  };
		\end{tikzpicture}
	\end{center}	
 
The following theorem states that there is a domain containing every natural number.

\begin{theorem}[Subdomain property of natural numbers, \ref{def:dom-prop}, \ref{def:ind-nat}]\ \\ \label{thm:subdomain-natural}
  $$ M(N)  $$
\end{theorem}
  \begin{center}
		\begin{tikzpicture}[every node/.style = {draw, shape=rectangle,anchor=north},growth parent anchor=south,level distance=.4 cm,sibling distance = 4.7cm]
		\node{$\ds{M}(N) $}
				child{[level distance=.4cm]node{$\lall \sigma =\sigma \rall \lexi N \xi \rexi \nonmem \xi \sigma $}
				child{[level distance=.1cm]node{$s=s $}
				child{[level distance=.4cm]node{$ Ss $}
				child{[level distance=.1cm]node{$ Nx $}
                                  child{[level distance=.4cm]node{$\nonmem x s $}
                                  child{[level distance=.4cm]node{$\lall S\phi \rall \mem x \phi $}
				child{node{$  \mem x s$}
				child{node{$  $}
edge from parent node[right, draw=none] {\color{red} 2}
                                }
                                    edge from parent node[left, draw=none] {\color{red} \ref*{rule:all}}
edge from parent node[right, draw=none] {\color{red} 0-3}
                              }
                                    edge from parent node[left, draw=none] {\color{red} \ref*{def:ind-nat}-\ref*{rule:subs}}
edge from parent node[right, draw=none] {\color{red} 1}
                                  }
                                }
                                    edge from parent node[left, draw=none] {\color{red} \ref*{rule:all}-\ref*{rule:exist}}
edge from parent node[right, draw=none] {\color{red} 1-2}
                          }
                        }
                                    edge from parent node[left, draw=none] {\color{red} \ref*{thm:succ-funct}-\ref*{rule:exist}}
edge from parent node[right, draw=none] {\color{red} 0}
                        }
                                    edge from parent node[left, draw=none] {\color{red} \ref*{def:dom-prop}-\ref*{rule:subs}}
edge from parent node[right, draw=none] {\color{red} 0}
                  };
		\end{tikzpicture}
	\end{center}	

        The next theorem gives an object with inductive domain containing precisely the natural numbers, floating
        with each natural number to the successor of the number.

\begin{theorem}[Minimal successor function, \ref{def:dom-prop}, \ref{def:id-ext-succ}, \ref{def:ref-eq}, \ref{def:ind-nat}]\ \\ \label{thm:successor-function}
$$\lexi T \sigma \rexi  {T}\sigma $$
\end{theorem}
\begin{center}
		\begin{tikzpicture}[every node/.style = {draw, shape=rectangle,anchor=north},growth parent anchor=south,level distance=.4 cm,sibling distance = 5.1cm]
		\node{$\lall T \sigma \rall  \ds{T}\sigma $}
				child{[level distance=.1cm]node{$ M(N)$}
				child{[level distance=.4cm]node{$\lexi N\xi \rexi \lall I \eta  \rall \ds{H}\eta \xi $}
				child{[level distance=.1cm]node{$Nx$}
				child{[level distance=.4cm]node{$\lall  I \eta \rall \ds{H}\eta x$}
                                  child{[level distance=.1cm]node{$ \mem x  f $}
                                  child{[level distance=.4cm]node{$\lall \mem \xi f \rall \lexi \mem \eta f \rexi H\eta \xi $}
                                  child{[level distance=.1cm]node{$  \mem y f $}
                                    child{[level distance=.4cm]node{$ Hyx $}
                                    child{node{$ Iy $}
                                    child{node{$ \ds{H}yx $}
                                    child{node{$  $}
edge from parent node[right, draw=none] {\color{red} 2}
                    }
                                    edge from parent node[left, draw=none] {\color{red} \ref*{rule:all}}
edge from parent node[right, draw=none] {\color{red} 0-5}
                    }
                                   edge from parent node[left, draw=none] {\color{red} \ref*{thm:abb-ref-eq}-\ref*{rule:contradiction}}
                    }
                    }
                                    edge from parent node[left, draw=none] {\color{red} \ref*{rule:all}-\ref*{rule:exist}}
edge from parent node[right, draw=none] {\color{red} 0-1}
                                                                      }
                                                                      }
                                    edge from parent node[left, draw=none] {\color{red} \ref*{thm:infinity}-\ref*{rule:contradiction}-\ref*{rule:exist}}
                                                                      }
                                                                      }
                                    edge from parent node[left, draw=none] {\color{red} \ref*{rule:exist}}
edge from parent node[right, draw=none] {\color{red} 0}
                                                                      }
                                                                      }
                                  child{[level distance=.4cm]node{$\lall N\xi \rall \lexi I \eta \rexi {H}\eta \xi $}
                                  child{[level distance=.1cm]node{$ D(N) s$}
                                  child{[level distance=.1cm, sibling distance=3.1cm]node{$\lall N\xi \rall \lexi I s(\xi) \rexi Hs(\xi)\xi $}
                                  child{[level distance=.4cm]node{${T}s$}
                                    child{node{$ \ds{T}s $}
                                  child{[level distance=.1cm]node{$ $}
edge from parent node[right, draw=none] {\color{red} 1}
                    }
                                    edge from parent node[left, draw=none] {\color{red} \ref*{rule:all}}
edge from parent node[right, draw=none] {\color{red} 0-5}
                    }
                    }
                                  child{[level distance=.4cm]node{$\ds{T}s$}
                                    child{node{$ \lall D(N) s \rall \lexi N\xi \rexi \ds{H}s (\xi)\xi $}
                                    child{[level distance=.1cm]node{$ Nx  $}
                                    child{[level distance=.4cm]node{$ \ds{H}s (x)x $}
                                    child{node{$ {H}s (x)x $}
                                  child{[level distance=.1cm]node{$ $}
edge from parent node[right, draw=none] {\color{red} 1}
                    }
                                    edge from parent node[left, draw=none] {\color{red} \ref*{rule:all}-\ref*{rule:exist}}
edge from parent node[right, draw=none] {\color{red} 1-4}
                    }
                    }
                                    edge from parent node[left, draw=none] {\color{red} \ref*{rule:all}-\ref*{rule:exist}}
edge from parent node[right, draw=none] {\color{red} 0-3}
                    }
                                    edge from parent node[left, draw=none] {\color{red} \ref*{def:ind-nat}-\ref*{rule:subs}}
edge from parent node[right, draw=none] {\color{red} 0}
                    }
                    }
                    }
                                    edge from parent node[left, draw=none] {\color{red} \ref*{thm:choice}-\ref*{rule:exist}}
edge from parent node[right, draw=none] {\color{red} 0-1}
                          }
                        }
     edge from parent node[left, draw=none] {\color{red} \ref*{thm:subdomain-natural}}
                  };
		\end{tikzpicture}
	\end{center}	

\begin{definition}[Successor function]\label{def:successor-function}
\end{definition}
\begin{center}
\begin{tikzpicture}[every node/.style = {draw, shape=rectangle,anchor=north},growth parent anchor=south,level distance=.4 cm,sibling distance = 4.7cm]]
  \node{$T \sigma \prop  T_{\ref*{def:ind-nat}}\sigma$}
  child{[level distance=.1cm]node{$Ts$}
       edge from parent node[left, draw=none] {\color{red} \ref*{thm:successor-function}}
  };
    \end{tikzpicture}
\end{center}

The object $s$ from Definition \ref{def:successor-function} has been constructed so that
it floats with every natural number  to the successor of the number.  The following theorem confirms that.

  \begin{theorem}[Auxiliary theorem, \ref{def:id-ext-succ}, \ref{def:ind-nat}, \ref{def:successor-function}]\ \\ \label{thm:nat-succ}
$$ \lall N\xi  \rall H{s(\xi)}\xi$$
\end{theorem}
\begin{center}
		\begin{tikzpicture}[every node/.style = {draw, shape=rectangle,anchor=north},growth parent anchor=south,level distance=.4 cm,sibling distance = 4.7cm]
		\node{$ \lexi N \xi  \rexi \ds{H}s(\xi)\xi$}
 			child{[level distance=.1cm]node{$ Nx $}
				child{[level distance=.4cm]node{$\ds{H} s(x) x $}
				child{[level distance=.4cm]node{$ Ts $}
				child{[level distance=.4cm]node{$ \lall N\xi \rall Hs(\xi)\xi $}
                                 child{[level distance=.4cm]node{$ Hs(x)x $}
                                  child{[level distance=.4cm]node{$  $}
edge from parent node[right, draw=none] {\color{red} 3}
                                  }
                                  edge from parent node[left, draw=none] {\color{red} \ref*{rule:exist}-\ref*{rule:all}}
                                  edge from parent node[right, draw=none] {\color{red} 0-3}
                                  }
                          edge from parent node[left, draw=none] {\color{red} \ref*{def:ind-nat}-\ref*{rule:subs}-\ref*{rule:exist}}
edge from parent node[right, draw=none] {\color{red} 0}
                          }
                          edge from parent node[left, draw=none] {\color{red} \ref*{def:successor-function}}
                         }
                    }
                                    edge from parent node[left, draw=none] {\color{red} \ref*{rule:exist}}
edge from parent node[right, draw=none] {\color{red} 0}
                  };
		\end{tikzpicture}
	\end{center}	

     There are five peano axioms. The first Peano axiom, Theorem \ref{thm:peanoi}, states that zero is a natural number.

  \begin{theorem}[Peano Axiom I, \ref{def:id-ext-succ}, \ref{def:empty-object}, \ref{def:ind-nat}]\ \\ \label{thm:peanoi}
$$ {N} o$$
\end{theorem}
\begin{center}
		\begin{tikzpicture}[every node/.style = {draw, shape=rectangle,anchor=north},growth parent anchor=south,level distance=.4 cm,sibling distance = 4.7cm]
		\node{$ \ds{N} o$}
 			child{[level distance=.4cm]node{$\lexi S\phi \rexi \nonmem o \phi $}
				child{[level distance=.1cm]node{$Sr $}
				child{[level distance=.4cm]node{$ \nonmem o r $}
				child{[level distance=.4cm]node{$ \lexi \mem o r \rexi \lall \mem \xi r\rall \lexi H\eta \xi \rexi \mem \eta r $}
                                  child{[level distance=.4cm]node{$\mem o r $}
                                  child{[level distance=.4cm]node{$  $}
edge from parent node[right, draw=none] {\color{red} 2}
                                  }
                                    edge from parent node[left, draw=none] {\color{red} \ref*{rule:exist}}
edge from parent node[right, draw=none] {\color{red} 0}
                                }
                                    edge from parent node[left, draw=none] {\color{red} \ref*{def:ind-nat}-\ref*{rule:subs}}
edge from parent node[right, draw=none] {\color{red} 1}
                          }
                         }
                                    edge from parent node[left, draw=none] {\color{red} \ref*{rule:exist}}
edge from parent node[right, draw=none] {\color{red} 0}
                    }
                                    edge from parent node[left, draw=none] {\color{red} \ref*{def:ind-nat}-\ref*{rule:subs}}
edge from parent node[right, draw=none] {\color{red} 0}
                  };
		\end{tikzpicture}
	\end{center}	

Induction is a proof method. If zero satisfies an attribute, and with every natural number
that satisfies the attribute also its successor satisfies the attribute, then every natural number satisfies this attribute. 
This is stated in the following theorem. It is proved by producing an auxiliary inductive function.

\begin{theorem}[Induction,
  \ref{def:dom-prop}, \ref{def:id-ext-succ}, \ref{def:ind-nat}, \ref{def:successor-function}]\ \\ \label{thm:induction}
$$\lall \Gamma 0\rall \lall \lall \Gamma \xi \rall N\xi \rall
\lall \lall \Gamma \xi \rall \Gamma s(\xi)\rall \lall N \eta  \rall {\Gamma}\eta$$
\end{theorem}
  \begin{center}
		\begin{tikzpicture}[every node/.style = {draw, shape=rectangle,anchor=north},growth parent anchor=south,level distance=.4 cm,sibling distance = 4.7cm]
                  \node{$\lexi G 0\rexi \lexi \lall G \xi \rall N\xi \rexi
                    \lexi \lall G \xi \rall G s(\xi)\rexi \lexi N \eta  \rexi \ds{G}\eta$}
 			child{[level distance=.1cm]node{$G o $}
			child{[level distance=.1cm]node{$\lall G \xi \rall N\xi  $}
				child{[level distance=.1cm]node{$ \lall G \xi \rall G s(\xi)$}
				child{[level distance=.1cm]node{$ N y$}
				child{[level distance=.4cm, sibling distance=4cm]node{$ \ds{G}y $}
				child{[level distance=.4cm]node{$M(G)$}
				child{[level distance=.4cm]node{$D(G)r  $}
                                  child{[level distance=.1cm]node{$\lall G \xi \rall \mem \xi r$}
				child{[level distance=.4cm]node{$ \lall \mem \xi r \rall G \xi $}
				child{[level distance=.1cm,sibling distance=5.2cm]node{$ \mem o r  $}
				child{[level distance=.4cm]node{$ \ds{S}r  $}
				child{node{$ \lall \mem o r \rall \lexi \mem \xi r \rexi \lall \mem \eta r \rall \ds{H}\eta \xi $}
                                  child{[level distance=.1cm]node{$  \mem x r    $}
                                    child{[level distance=.4cm]node{$ \lall \mem \eta r \rall \ds{H}\eta x   $}
				child{node{$ G x  $}
				child{node{$ G s(x)  $}
				child{node{$ \mem {s(x)} r  $}
				child{node{$ Nx  $}
				child{node{$ Hs(x)x  $}
				child{node{$ \ds{H}s(x)x  $}
				child{node{$  $}
edge from parent node[right, draw=none] {\color{red} 1}
                                }
                                    edge from parent node[left, draw=none] {\color{red} \ref*{rule:all}}
edge from parent node[right, draw=none] {\color{red} 2-5}
                                }
                                    edge from parent node[left, draw=none] {\color{red} \ref*{thm:nat-succ}-\ref*{rule:all}}
edge from parent node[right, draw=none] {\color{red} 0}
                                }
                                    edge from parent node[left, draw=none] {\color{red} \ref*{rule:all}}
edge from parent node[right, draw=none] {\color{red} 2-15}
                                }
                                    edge from parent node[left, draw=none] {\color{red} \ref*{rule:all}}
edge from parent node[right, draw=none] {\color{red} 0-8}
                                }
                                    edge from parent node[left, draw=none] {\color{red} \ref*{rule:all}}
edge from parent node[right, draw=none] {\color{red} 0-12}
                                }
                                    edge from parent node[left, draw=none] {\color{red} \ref*{rule:all}}
edge from parent node[right, draw=none] {\color{red} 1-5}
                                }
                                }
                                    edge from parent node[left, draw=none] {\color{red} \ref*{rule:all}-\ref*{rule:exist}}
edge from parent node[right, draw=none] {\color{red} 0-2}
                                }
                               edge from parent node[left, draw=none] {\color{red} \ref*{def:ind-nat}- \ref*{rule:subs}}
edge from parent node[right, draw=none] {\color{red} 0}
                                }
}
                                  child{[level distance=.4cm]node{$ {S}r $}
				child{[level distance=.4cm]node{$ \lall S\phi \rall  \mem y \phi $}
				child{[level distance=.4cm]node{$ \mem y r  $}
				child{[level distance=.4cm]node{$ G y  $}
				child{[level distance=.4cm]node{$  $}
edge from parent node[right, draw=none] {\color{red} 9}
                                }
                                 edge from parent node[left, draw=none] {\color{red} \ref*{rule:exist}}
edge from parent node[right, draw=none] {\color{red} 0-4}
                                }
                                    edge from parent node[left, draw=none] {\color{red} \ref*{rule:all}}
edge from parent node[right, draw=none] {\color{red} 0-1}
                                }
                         edge from parent node[left, draw=none] {\color{red} \ref*{def:ind-nat}-\ref*{rule:subs}}
edge from parent node[right, draw=none] {\color{red} 7}
                                }
                                }
                                    edge from parent node[left, draw=none] {\color{red} \ref*{rule:all}}
edge from parent node[right, draw=none] {\color{red} 1-8}
                              }
                                }
                                    edge from parent node[left, draw=none] {\color{red} \ref*{def:dom-prop}-\ref*{rule:subs}-\ref*{rule:exist}}
edge from parent node[right, draw=none] {\color{red} 0}
                                }
                                    edge from parent node[left, draw=none] {\color{red} \ref*{thm:exist-identity}-\ref*{rule:all}-\ref*{rule:exist}}
edge from parent node[right, draw=none] {\color{red} 0}
                              }
                                    edge from parent node[left, draw=none] {\color{red} \ref*{thm:nested-att}-\ref*{thm:subdomain-natural}-\ref*{rule:all}}
edge from parent node[right, draw=none] {\color{red} 3}
                                  }
                                }
                        }
                    }
                    }
                                    edge from parent node[left, draw=none] {\color{red} \ref*{rule:exist}}
edge from parent node[right, draw=none] {\color{red} 0}
                  };
		\end{tikzpicture}
	\end{center}	

The next theorem shows  a uniqueness result for the successor. It is an application of the rule of function uniqueness.

\begin{theorem}  [Uniqueness of successor, \ref{def:id-ext-succ}]\ \\ 	\label{thm:successor-uniqueness}
$$\lall \phi=\phi\rall  \lall H\psi \phi\rall \lall H\rho \phi \rall \psi= \rho$$
\end{theorem}
\begin{center}
		\begin{tikzpicture}[every node/.style = {draw, shape=rectangle,anchor=north},growth parent anchor=south,level distance=.4 cm,sibling distance = 4.7cm]
                  \node{$\lexi \phi=\phi\rexi  \lexi H\psi \phi\rexi \lexi H\rho \phi \rexi \psi\ds{\eq} \rho
                    $}
		child{[level distance=.1cm]node{$f=f$} 
                  child{[level distance=.1cm]node{$Hgf
                      $} 
					child{[level distance=.1cm]node{$Hhf$} 
					child{[level distance=.4cm]node{$g\ds{\eq} h$} 
                                          child{[sibling distance = 4.2cm,level distance=.1cm]node{$\lall Ef\xi \rall \mem \xi g$}
                                              child{[level distance=.1cm]node{$\lall \mem \xi g \rall Ef\xi $}
                                              child{[level distance=.4cm]node{$ Cg$}
                                          child{[sibling distance = 6.3cm,level distance=.1cm]node{$ \lall Ef\xi \rall \mem \xi h$}
                                              child{[level distance=.1cm]node{$\lall \mem \xi h \rall Ef\xi $}
                                              child{[level distance=.4cm]node{$ Ch$}
										child{[level distance=.1cm]node{$g(x)\ds{\eq} h(x)$} 
                                         child{[sibling distance=2.5cm,level distance=.1cm]node{$\lall \nonmem x g\rall \mem x h$} 
					child{[level distance=.4cm]node{$\nonmem x g$}
                                          child{[level distance=0.4cm]
                                            node{$\mem x h $}
                                            child{[level distance=.4cm]node{$  Efx  $}
                                            child{[level distance=.4cm]node{$  \mem x g $}
                                              child{[level distance=.4cm]node{$ $}
edge from parent node[right, draw=none] {\color{red} 3}
                                              }
edge from parent node[left, draw=none] {\color{red} \ref*{rule:all}}
edge from parent node[right, draw=none] {\color{red} 0-10}
}
edge from parent node[left, draw=none] {\color{red} \ref*{rule:all}}
edge from parent node[right, draw=none] {\color{red} 0-5}
}
edge from parent node[left, draw=none] {\color{red} \ref*{rule:all}}
edge from parent node[right, draw=none] {\color{red} 0-1}
}}
                                            child{[level distance=.4cm]node{$  \mem x g $}
                                              child{[level distance=.4cm]node{$ Efx$}
                                              child{[level distance=.4cm]node{$ \mem x h$}
                                              child{[level distance=.4cm]node{$ g(x)=x$}
                                              child{[level distance=.4cm]node{$ h(x)=x$}
                                              child{[level distance=.4cm]node{$ g(x)=h(x)$}
                                              child{[level distance=.4cm]node{$ $}
edge from parent node[right, draw=none] {\color{red} 7}
                                              }
edge from parent node[left, draw=none] {\color{red} \ref*{rule:subf}}
edge from parent node[right, draw=none] {\color{red} 0-1}
                                            }
edge from parent node[left, draw=none] {\color{red}  \ref*{def:id-ext-succ}-\ref*{rule:subs}-\ref*{rule:all}}
edge from parent node[right, draw=none] {\color{red} 1-6}
                                          }
edge from parent node[left, draw=none] {\color{red} \ref*{def:id-ext-succ}-\ref*{rule:subs}-\ref*{rule:all}}
edge from parent node[right, draw=none] {\color{red} 2-8}
                                        }
edge from parent node[left, draw=none] {\color{red} \ref*{rule:all}}
edge from parent node[right, draw=none] {\color{red} 0-6}
                                  }
edge from parent node[left, draw=none] {\color{red} \ref*{rule:all}}
edge from parent node[right, draw=none] {\color{red} 0-7}
                              }
}
}
edge from parent node[left, draw=none] {\color{red} \ref*{thm:function-uniqueness}-\ref*{rule:contradiction}-\ref*{rule:all}-\ref*{rule:exist}}
edge from parent node[right, draw=none] {\color{red} 6}
}
}
                                }
edge from parent node[left, draw=none] {\color{red} \ref*{def:id-ext-succ}-\ref*{rule:subs}-\ref*{rule:exist}}
edge from parent node[right, draw=none] {\color{red} 4}
}
}
                                              }
edge from parent node[left, draw=none] {\color{red} \ref*{def:id-ext-succ}-\ref*{rule:subs}-\ref*{rule:exist}}
edge from parent node[right, draw=none] {\color{red} 2}
}
}
}
}
edge from parent node[left, draw=none] {\color{red} \ref*{rule:exist}}
edge from parent node[right, draw=none] {\color{red} 0}
};
		\end{tikzpicture}
	\end{center}	

The second Peano axiom, Theorem \ref{thm:peanoii},  states that the successor of a natural number is
again a natural number. This is inherited from the general inductive function, since natural numbers are precisely
the objects in the intersection of all inductive domains.

  \begin{theorem}[Peano Axiom II, \ref{def:id-ext-succ}, \ref{def:ind-nat}, \ref{def:successor-function}]\ \\ \label{thm:peanoii}
  $$ \lall N \xi \rall {N}s(\xi)$$
   \end{theorem}
	\begin{center}
		\begin{tikzpicture}[every node/.style = {draw, shape=rectangle,anchor=north},growth parent anchor=south,level distance=.4 cm,sibling distance = 4.7cm]
		\node{$ \lexi N \xi \rexi \ds{N}s(\xi) $}
                                  child{[level distance=.1cm]node{$ Nx  $}
                                  child{[level distance=.4cm]node{$ \ds{N}s(x)$}
                                  child{[level distance=.4cm]node{$\lall S\phi \rall \mem x \phi $}
 			child{node{$\lexi S\phi \rexi \nonmem {s(x)} \phi $}
				child{[level distance=.1cm]node{$Sr $}
				child{[level distance=.4cm]node{$ \nonmem {s(x)} r $}
                                  child{[level distance=.4cm]node{$\mem x r $}
				child{[level distance=.4cm]node{$ \lexi \mem o r \rexi \lall \mem \xi r\rall \lexi \mem \eta r \rexi  {H}\eta \xi$}
                                   child{[level distance=.1cm]node{$ \mem y r $}
                                  child{[level distance=.4cm]node{$ Hyx $}
                              child{[level distance=.4cm]node{$Hs(x)x $}
                                  child{[level distance=.4cm]node{$ y=s(x)$}
                                  child{[level distance=.4cm]node{$ \nonmem y r $}
                                  child{[level distance=.4cm]node{$  $}
edge from parent node[right, draw=none] {\color{red} 4}
                                  }
                                    edge from parent node[left, draw=none] {\color{red} \ref*{rule:subf}}
edge from parent node[right, draw=none] {\color{red} 0-6}
                                  }
                                    edge from parent node[left, draw=none] {\color{red} \ref*{thm:successor-uniqueness}-\ref*{rule:contradiction}-\ref*{rule:all}}
edge from parent node[right, draw=none] {\color{red} 0-1}
                                  }
                                    edge from parent node[left, draw=none] {\color{red} \ref*{thm:nat-succ}}
edge from parent node[right, draw=none] {\color{red} 9}
                                  }
                                  }
                                    edge from parent node[left, draw=none] {\color{red} \ref*{rule:exist}-\ref*{rule:all}}
edge from parent node[right, draw=none] {\color{red} 0-1}
                                  }
                                    edge from parent node[left, draw=none] {\color{red} \ref*{def:ind-nat}-\ref*{rule:subs}}
edge from parent node[right, draw=none] {\color{red} 2}
                                  }
                                    edge from parent node[left, draw=none] {\color{red} \ref*{rule:all}}
edge from parent node[right, draw=none] {\color{red} 1-3}
                                  }
                                  }
                                    edge from parent node[left, draw=none] {\color{red} \ref*{rule:exist}}
edge from parent node[right, draw=none] {\color{red} 0}
                                }
                                    edge from parent node[left, draw=none] {\color{red}\ref*{def:ind-nat}-\ref*{rule:subs}}
edge from parent node[right, draw=none] {\color{red} 1}
                          }
                                    edge from parent node[left, draw=none] {\color{red} \ref*{def:ind-nat}-\ref*{rule:subs}}
edge from parent node[right, draw=none] {\color{red} 1}
                          }
                    }
                                    edge from parent node[left, draw=none] {\color{red} \ref*{rule:exist}}
edge from parent node[right, draw=none] {\color{red} 0}
                  };
		\end{tikzpicture}
	\end{center}	


We postpone the discussion of the Peano axiom that is traditionally called the third.
The fourth Peano axiom, Theorem \ref{thm:peanoiv}, states that zero is not the successor of any natural number. Its proof
uses the following theorem, that states that the successor of a natural number contains at least this natural number.
The successor therefore can not have an empty domain.
 
  \begin{theorem}[Natural number in successor, \ref{def:id-ext-succ}, \ref{def:ind-nat}, \ref{def:successor-function}]\ \\ \label{thm:nat-num-in-suc}
  $$ \lall N\xi \rall \mem \xi {s(\xi)}$$
\end{theorem}

\begin{center}
		\begin{tikzpicture}[every node/.style = {draw, shape=rectangle,anchor=north},growth parent anchor=south,level distance=.4 cm,sibling distance = 4.7cm]
		\node{$ \lexi N\xi \rexi \nonmem x {s(x)}$}
 			child{[level distance=.1cm]node{$ Nx $}
				child{[level distance=.4cm]node{$ \nonmem x {s(x)} $}
                                child{[level distance=.4cm]node{$ Hs(x) x  $}
                                  child{[level distance=.1cm]node{$ \lall Ex\xi \rall  \mem \xi {s(x)} $}
                                  child{[level distance=.4cm]node{$ Exx $}
                                  child{[level distance=.4cm]node{$ \mem x {s(x)} $}
                                  child{[level distance=.4cm]node{$ $}
edge from parent node[right, draw=none] {\color{red} 4}
                                  }
                                    edge from parent node[left, draw=none] {\color{red} \ref*{rule:all}}
edge from parent node[right, draw=none] {\color{red} 0-1}
                                  }
                                  }
                                 child{[level distance=.4cm]node{$ \ds{E}xx $}
                                  child{[level distance=.4cm]node{$ \lexi \nonmem x x \rexi x\ds{\eq} x $}
                                  child{[level distance=.1cm]node{$x\ds{\eq} x $}
                                  child{[level distance=.4cm]node{$  $}
                                  }
                                    edge from parent node[left, draw=none] {\color{red} \ref*{rule:exist}}
edge from parent node[right, draw=none] {\color{red} 0}
                                  }
                                  edge from parent node[left, draw=none] {\color{red}
                                     \ref*{def:id-ext-succ}-\ref*{rule:subs}}
edge from parent node[right, draw=none] {\color{red} 0}
                                  }
                                }
                                   edge from parent node[left, draw=none] {\color{red} \ref*{def:id-ext-succ}-\ref*{rule:subs}-\ref*{rule:exist}}
edge from parent node[right, draw=none] {\color{red} 0}
                                  }
                                    edge from parent node[left, draw=none] {\color{red} \ref*{thm:nat-succ}}
edge from parent node[right, draw=none] {\color{red} 1}
                                  }
                    }
                                    edge from parent node[left, draw=none] {\color{red} \ref*{rule:exist}}
edge from parent node[right, draw=none] {\color{red} 0}
                  };
		\end{tikzpicture}
	\end{center}	

  \begin{theorem}[Peano Axiom IV, \ref{def:empty-object}, \ref{def:ind-nat}, \ref{def:successor-function}]\ \\ \label{thm:peanoiv}
$$ \lall N \xi \rall s(\xi)\ds{\eq}0 $$
\end{theorem}
	\begin{center}
		\begin{tikzpicture}[every node/.style = {draw, shape=rectangle,anchor=north},growth parent anchor=south,level distance=.4 cm,sibling distance = 4.7cm]
		\node{$\lexi N \xi \rexi s(\xi){\eq}0$}
                                  child{[level distance=.1cm]node{$ Nx  $}
                                  child{[level distance=.4cm]node{$ s(x)=o$}
                                 child{[level distance=.4cm]node{$\mem x {s(x)} $}
				child{[level distance=.4cm]node{$ \mem x o  $}
                                  child{[level distance=.4cm]node{$\ 0 o $}
                                  child{[level distance=.4cm]node{$ \nonmem x o $}
                                  child{[level distance=.4cm]node{$ $}
edge from parent node[right, draw=none] {\color{red} 2}
                                  }
                                    edge from parent node[left, draw=none] {\color{red} \ref*{def:empty}-\ref*{def:empty-object}-\ref*{rule:subs}-\ref*{rule:contradiction}-\ref*{rule:all}}
edge from parent node[right, draw=none] {\color{red} 0}
                                  }
                                    edge from parent node[left, draw=none] {\color{red} \ref*{def:empty-object}}
                                  }
                                    edge from parent node[left, draw=none] {\color{red} \ref*{rule:subf}}
edge from parent node[right, draw=none] {\color{red} 0-1}
                                  }
                                    edge from parent node[left, draw=none] {\color{red} \ref*{thm:nat-num-in-suc}-\ref*{rule:all}}
edge from parent node[right, draw=none] {\color{red} 1}
                                  }
                    }
                                    edge from parent node[left, draw=none] {\color{red} \ref*{rule:exist}}
edge from parent node[right, draw=none] {\color{red} 0}
                  };
		\end{tikzpicture}
	\end{center}	

The fifth Peano axiom, Theorem \ref{thm:peanov}, states the inductive principle. It is a variant of
Theorem \ref{thm:induction}, which is prominently used in the proof.

\begin{theorem}[Peano Axiom V, \ref{def:ind-nat}, \ref{def:successor-function}]\ \\ \label{thm:peanov}
  $$ \lall \Sigma o \rall \lall \lall N\xi \rall \lall \Sigma \xi \rall \Sigma s(\xi)\rall \lall N\xi \rall
                    {\Sigma} \xi $$
  \end{theorem}

	\begin{center}
		\begin{tikzpicture}[every node/.style = {draw, shape=rectangle,anchor=north},growth parent anchor=south,level distance=.4 cm,sibling distance = 5.8cm]
                  \node{$ \lexi G o \rexi \lexi \lall N\xi \rall \lall G \xi \rall G s(\xi)\rexi \lexi N\xi \rexi \ds{G} \xi $}
                                  child{[level distance=.1cm]node{$ G o  $}
                                  child{[level distance=.1cm]node{$ \lall N\xi \rall \lall G \xi \rall G s(\xi))$}
                                  child{[level distance=.1cm]node{$Nx $}
 			child{[level distance=.4cm]node{$\ds{G} x $}
				child{[level distance=.4cm]node{$A\xi : \lexi N\xi \rexi G \xi $}
                                  child{[level distance=.1cm, sibling distance=4.3cm]node{
                                      $\lall A 0\rall \lall \lall A \xi \rall N\xi \rall
                                      \lall \lall A \xi \rall A s(\xi)\rall \lall N \eta  \rall A \eta $
                                    }
                                child{[level distance=.4cm]node{$ \lexi A \xi \rexi \ds{A} s(\xi)$}
                                  child{[level distance=.1cm]node{$ Ax  $}
                                  child{[level distance=.4cm]node{$ \ds{A} s(x)$}
                                  child{[level distance=.4cm]node{$ \lexi Nx\rexi G x$}
                                  child{[level distance=.1cm]node{$ Nx $}
                                  child{[level distance=.4cm]node{$ G x $}
                                  child{[level distance=.4cm]node{$ N s(x) $}
                                  child{[level distance=.4cm]node{$ G s(x) $}
                                  child{[level distance=.4cm]node{$ \lall N s(x)\rall \ds{G}s(x) $}
                                  child{[level distance=.4cm]node{$  \ds{G}s(x) $}
                                  child{[level distance=.4cm]node{$  $}
edge from parent node[right, draw=none] {\color{red} 2}
                                  }
                                    edge from parent node[left, draw=none] {\color{red} \ref*{rule:all}}
edge from parent node[right, draw=none] {\color{red} 0-2}
                                  }
                                    edge from parent node[left, draw=none] {\color{red} \ref*{rule:subs}}
edge from parent node[right, draw=none] {\color{red} 5-9}
                                  }
                                  edge from parent node[left, draw=none] {\color{red} \ref*{rule:all}}
edge from parent node[right, draw=none] {\color{red} 1-2-11}
                                  }
                                    edge from parent node[left, draw=none] {\color{red}  \ref*{thm:peanoii}-\ref*{rule:all}}
edge from parent node[right, draw=none] {\color{red} 1}
                                  }
                                  }
                                    edge from parent node[left, draw=none] {\color{red} \ref*{rule:exist}}
edge from parent node[right, draw=none] {\color{red} 0}
                                  }
                                    edge from parent node[left, draw=none] {\color{red} \ref*{rule:subs}}
edge from parent node[right, draw=none] {\color{red} 1-4}
                                  }
                                  }
                                    edge from parent node[left, draw=none] {\color{red} \ref*{rule:exist}}
edge from parent node[right, draw=none] {\color{red} 0}
                                  }
                                  }
                                  child{[level distance=.1cm,sibling distance=3.7cm]node{$ \lall A\xi \rall As(\xi)$}
                                  child{[level distance=.4cm]node{$ \ds{A} o$}
                                  child{[level distance=.4cm]node{$ \lall No \rall \ds{G} o  $}
                                  child{[level distance=.4cm]node{$ N o$}
                                  child{[level distance=.4cm]node{$ \ds{G} o$}
                                  child{[level distance=.4cm]node{$  $}
edge from parent node[right, draw=none] {\color{red} 10}
                                  }
                                    edge from parent node[left, draw=none] {\color{red} \ref*{rule:all}}
edge from parent node[right, draw=none] {\color{red} 0-1}
                                  }
                                    edge from parent node[left, draw=none] {\color{red} \ref*{thm:peanoi}}
                                  }
                                    edge from parent node[left, draw=none] {\color{red} \ref*{rule:subs}}
edge from parent node[right, draw=none] {\color{red} 0-3}
                                  }
                                  }
                                  child{[level distance=.1cm, sibling distance=2.3cm]node{$A o $}
       				child{[level distance=.4cm]node{$ \lall A\xi \rall N\xi$}
				                                 child{[level distance=.4cm]node{$ A x $}
                                  child{[level distance=.4cm]node{$ \lexi Nx \rexi G x $}
                                  child{[level distance=.4cm]node{$ G x $}
                                  child{[level distance=.4cm]node{$ $}
edge from parent node[right, draw=none] {\color{red} 8}
                                  }
                                    edge from parent node[left, draw=none] {\color{red} \ref*{rule:exist}}
edge from parent node[right, draw=none] {\color{red} 0}
                                  }
                                    edge from parent node[left, draw=none] {\color{red} \ref*{rule:subs}}
edge from parent node[right, draw=none] {\color{red} 0-5}
                                  }
                                    edge from parent node[left, draw=none] {\color{red} \ref*{rule:all}}
edge from parent node[right, draw=none] {\color{red} 0-1-2-3-6}
                                  }
                                  }
                                  child{[level distance=.4cm]node{$ \lexi A \xi \rexi \ds{N} \xi $}
                                  child{[level distance=.1cm]node{$ Ax  $}
                                  child{[level distance=.4cm]node{$ \ds{N}x$}
                                  child{[level distance=.4cm]node{$ \lexi Nx\rexi G x$}
                                  child{[level distance=.4cm]node{$ Nx $}
                          child{[level distance=.4cm]node{$ $}
edge from parent node[right, draw=none] {\color{red} 2}
                                  }
                                    edge from parent node[left, draw=none] {\color{red} \ref*{rule:exist}}
edge from parent node[right, draw=none] {\color{red} 0}
                                  }
                                    edge from parent node[left, draw=none] {\color{red} \ref*{rule:subs}}
edge from parent node[right, draw=none] {\color{red} 1-6}
                                  }
                                  }
                                    edge from parent node[left, draw=none] {\color{red} \ref*{rule:exist}}
edge from parent node[right, draw=none] {\color{red} 0}
                                  }
                                  }                                                              
                                  }
                                  }                               
                                   edge from parent node[left, draw=none] {\color{red} \ref*{thm:induction}}
                                  }
                                    edge from parent node[left, draw=none] {\color{red} \ref*{rule:con-pro}}
                                }
                          }
                          }
                    }
                                    edge from parent node[left, draw=none] {\color{red} \ref*{rule:exist}}
edge from parent node[right, draw=none] {\color{red} 0}
}
;
		\end{tikzpicture}
	\end{center}	

To prove the third Peano axiom we need some preparation. 
Theorem \ref{thm:predecessor} proves by induction that every natural number other than zero is the successor
of some natural number. The theorem is then used to prove Theorem \ref{thm:nat-num-id-f}, which
shows that every natural number is an identity object.


  \begin{theorem}[Predecessor, \ref{def:empty-object}, \ref{def:ind-nat}, \ref{def:successor-function}]\ \\ \label{thm:predecessor}
$$ \lall N \xi \rall \lall \xi\ds{\eq} o\rall \lexi N\eta \rexi s(\eta)\eq \xi $$
      \end{theorem}

    {\normalsize
    \begin{center}
		\begin{tikzpicture}[every node/.style = {draw, shape=rectangle,anchor=north},growth parent anchor=south,level distance=.4 cm,sibling distance = 4.7cm]
		\node{$ \lexi N \xi \rexi \lexi \xi\ds{\eq} o\rexi \lall N\eta \rall s(\eta)\ds{\eq} \xi $}
 			child{[level distance=.1cm]node{$ Nx $}
				child{[level distance=.1cm]node{$x\ds{\eq} o $}
				child{[level distance=.4cm]node{$ \lall N \eta  \rall s(\eta)\ds{\eq} x $}
                                  child{[level distance=.4cm]node{$ A\xi\prop  \lall \xi \ds{\eq} o\rall
                                      \lexi N\eta \rexi s(\eta)=\xi$}
                                    child{[level distance=.1cm]node{$  \lall Ao\rall \lall \lall N\xi \rall
                                        \lall A\xi \rall As(\xi)\rall \lall N\xi \rall A\xi$}
                                  child{[level distance=.4cm]node{$
                                      \lexi N\xi \rexi \lexi A\xi \rexi \ds{A}s(\xi)$}
                                  child{[level distance=.1cm]node{$ Ny$}
                                  child{[level distance=.4cm]node{$ \ds{A}s(y) $}
                                  child{[level distance=.4cm]node{$ \lexi s(y)\ds{\eq} 0\rexi \lall N\eta \rall s(\eta)\ds{\eq} s(y) $}
                                  child{[level distance=.1cm]node{$ s(y)\ds{\eq} s(y) $}
                                  child{[level distance=.4cm]node{$  $}
                                  }
                                    edge from parent node[left, draw=none] {\color{red} \ref*{rule:exist}-\ref*{rule:all}}
edge from parent node[right, draw=none] {\color{red} 0-2}
                                  }
                                    edge from parent node[left, draw=none] {\color{red} \ref*{rule:subs}}
edge from parent node[right, draw=none] {\color{red} 0-4}
                                  }
                                  }
                                    edge from parent node[left, draw=none] {\color{red} \ref*{rule:exist}}
edge from parent node[right, draw=none] {\color{red} 0}
                                  }
                                  }
                                  child{[level distance=.1cm, sibling distance=2.7cm]node{$ \lall N\xi \rall \lall A\xi \rall As(\xi) $}
                                  child{[level distance=.4cm]node{$  \ds{A} o$}
                                  child{[level distance=.1cm]node{$ o\ds{\eq} o$}
                                  child{[level distance=.4cm]node{$  $}
                                  }
                                    edge from parent node[left, draw=none] {\color{red} \ref*{rule:subs}-\ref*{rule:exist}}
edge from parent node[right, draw=none] {\color{red} 0-3}
                                  }
                                  }
                                  child{[level distance=.4cm]node{$ Ao $}
                                  child{[level distance=.4cm]node{$ Ax $}
                                  child{[level distance=.4cm]node{$  \lall x \ds{\eq} o\rall
                                      \lexi N\eta \rexi s(\eta)=x$}
                                  child{[level distance=.4cm]node{$ \lexi N\eta \rexi s(\eta)=x$}
                                  child{[level distance=.4cm]node{$  $}
edge from parent node[right, draw=none] {\color{red} 7}
                                  }
                                    edge from parent node[left, draw=none] {\color{red} \ref*{rule:all}}
edge from parent node[right, draw=none] {\color{red} 0-7}
                    }
                                    edge from parent node[left, draw=none] {\color{red} \ref*{rule:subs}}
edge from parent node[right, draw=none] {\color{red} 0-4}
                                  }
                                    edge from parent node[left, draw=none] {\color{red} \ref*{rule:all}}
edge from parent node[right, draw=none] {\color{red} 0-1-2-6}
                                  }
                                  }
                                  }
                                    edge from parent node[left, draw=none] {\color{red} \ref*{thm:peanov}}
                                  }
                                    edge from parent node[left, draw=none] {\color{red} \ref*{rule:con-pro}}
                                }
                          }
                        }
                                    edge from parent node[left, draw=none] {\color{red} \ref*{rule:exist}}
edge from parent node[right, draw=none] {\color{red} 0}
                  };
		\end{tikzpicture}
	\end{center}	
}


  \begin{theorem}[Natural number as identity object, \ref{def:id-ext-succ}, \ref{def:empty-object}, \ref{def:ind-nat}, \ref{def:successor-function}]\ \\ \label{thm:nat-num-id-f}
    $$\lall N \xi \rall C\xi $$
  \end{theorem}
    {
      \begin{center}
		\begin{tikzpicture}[every node/.style = {draw, shape=rectangle,anchor=north},growth parent anchor=south,level distance=.4 cm,sibling distance = 4.7cm]
		\node{$ \lexi N \xi \rexi \ds{C} \xi $}
 			child{[level distance=.1cm]node{$ Nx $}
				child{[level distance=.1cm]node{$\ds{C}x $}
                                  child{[level distance=.4cm]node{$x=o$}
                                  child{[level distance=.4cm]node{$Co$}
                                  child{[level distance=.4cm]node{$Cx $}
                              child{[level distance=.4cm]node{$  $}
edge from parent node[right, draw=none] {\color{red} 3}
                                  }
                                    edge from parent node[left, draw=none] {\color{red} \ref*{rule:subf}}
edge from parent node[right, draw=none] {\color{red} 0-1}
                                  }
                                    edge from parent node[left, draw=none] {\color{red} \ref*{thm:empty-id}}
                                  }
                                  }
                                  child{[level distance=.4cm, sibling distance=2.7cm]node{$ x\ds{\eq} o$}
                                 child{[level distance=.1cm]node{$ Ny $}
                                  child{[level distance=.4cm]node{$  s(y)=x$}
                               child{[level distance=.4cm]node{$ Hs(y)y $}
                                  child{[level distance=.4cm]node{$ H x y $}
                                  child{[level distance=.4cm]node{$ Cx$}
                                  child{[level distance=.4cm]node{$ $}
edge from parent node[right, draw=none] {\color{red} 6}
                                  }
                                    edge from parent node[left, draw=none] {\color{red} \ref*{def:id-ext-succ}-\ref*{rule:subs}-\ref*{rule:exist}}
edge from parent node[right, draw=none] {\color{red} 0}
                                  }
                                    edge from parent node[left, draw=none] {\color{red} \ref*{rule:subf}}
edge from parent node[right, draw=none] {\color{red} 0-1}
                                  }
                                    edge from parent node[left, draw=none] {\color{red} \ref*{thm:nat-succ}-\ref*{rule:all}}
edge from parent node[right, draw=none] {\color{red} 1}
                                  }
                                  }
                                    edge from parent node[left, draw=none] {\color{red}  \ref*{thm:predecessor}-\ref*{rule:all}-\ref*{rule:exist}}
edge from parent node[right, draw=none] {\color{red} 0-2}
                                  }
                                }
                        }
                                    edge from parent node[left, draw=none] {\color{red} \ref*{rule:exist}}
edge from parent node[right, draw=none] {\color{red} 0}
                  };
		\end{tikzpicture}
	\end{center}	
}

The following three theorems together with Theorem  \ref{thm:nat-num-in-suc} compare the domain of a natural number with
that of its successor.
While the natural number does not contain itself, Theorem \ref{thm:irreflexive-containment}, its successor
does contain the natural number, Theorem \ref{thm:nat-num-in-suc}. 
The domain of the successor is larger than the domain
of the natural number, Theorem \ref{thm:growth-nat-num}, but only by the natural number itself, Theorem \ref{thm:upper-successor}.

\begin{theorem}[Irreflexivity of containment, \ref{def:id-ext-succ}, \ref{def:ind-nat}]\ \\ \label{thm:irreflexive-containment}
  $$\lall N\xi \rall \nonmem \xi \xi $$
  \end{theorem}
\begin{center}
		\begin{tikzpicture}[every node/.style = {draw, shape=rectangle,anchor=north},growth parent anchor=south,level distance=.4 cm,sibling distance = 4.7cm]
		\node{$ \lexi N\xi \rexi \mem \xi \xi $}
 			child{[level distance=.1cm]node{$ Nx$}
				child{[level distance=.4cm]node{$\mem x x  $}
				child{[level distance=.4cm]node{$  C x  $}
				child{[level distance=.4cm]node{$ \nonmem x x  $}
				child{[level distance=.4cm]node{$ $}
edge from parent node[right, draw=none] {\color{red} 2}
                          }
                                    edge from parent node[left, draw=none] {\color{red} \ref*{thm:non-self-id}}
edge from parent node[right, draw=none] {\color{red} 0}
                    }                          
                              edge from parent node[left, draw=none] {\color{red} \ref*{thm:nat-num-id-f}}
edge from parent node[right, draw=none] {\color{red} 1}
                    }
                         }
                                    edge from parent node[left, draw=none] {\color{red} \ref*{rule:exist}}
edge from parent node[right, draw=none] {\color{red} 0}
                  };
		\end{tikzpicture}
	\end{center}	

  \begin{theorem}[Growth of natural number domain, \ref{def:ind-nat}, \ref{def:successor-function}]\ \\ \label{thm:growth-nat-num}
 $$\lall N\xi \rall \lall \mem \eta \xi \rall \mem \eta {s(\xi)}$$
\end{theorem}
	\begin{center}
		\begin{tikzpicture}[every node/.style = {draw, shape=rectangle,anchor=north},growth parent anchor=south,level distance=.4 cm,sibling distance = 4.7cm]
		\node{$ \lexi N\xi \rexi \lexi \mem \eta \xi \rexi \nonmem \eta {s(\xi)}$}
 			child{[level distance=.1cm]node{$Nx $}
				child{[level distance=.1cm]node{$\mem y x $}
				child{[level distance=.4cm]node{$ \nonmem y {s(x)} $}
                                 child{[level distance=.4cm]node{$ Hs(x)x $}
                                  child{[level distance=.1cm]node{$ \lall E x \xi \rall \mem \xi {s(x)} $}
                                  child{[level distance=.4cm]node{$ \ds{E}xy$}
                                  child{[level distance=.4cm]node{$ \lexi \nonmem y x \rexi x\ds{\eq} y $}
                                  child{[level distance=.4cm]node{$ \nonmem y x $}
                                  child{[level distance=.4cm]node{$ $}
edge from parent node[right, draw=none] {\color{red} 6}
                                  }
                                   edge from parent node[left, draw=none] {\color{red} \ref*{rule:exist}}
edge from parent node[right, draw=none] {\color{red} 0}
                                  }
                                    edge from parent node[left, draw=none] {\color{red} \ref*{rule:subs}}
edge from parent node[right, draw=none] {\color{red} 0}
                                  }
                                  }
                                  child{[level distance=.4cm]node{$ {E}xy$}
                                  child{[level distance=.4cm]node{$ \mem y {s(x)} $}
                                  child{[level distance=.4cm]node{$ $}
edge from parent node[right, draw=none] {\color{red} 4}
                                  }
                                    edge from parent node[left, draw=none] {\color{red} \ref*{rule:all}}
edge from parent node[right, draw=none] {\color{red} 0-1}
                                  }
                                  }
                                    edge from parent node[left, draw=none] {\color{red} \ref*{def:id-ext-succ}-\ref*{rule:subs}-\ref*{rule:exist}}
edge from parent node[right, draw=none] {\color{red} 0}
                                  }
                                    edge from parent node[left, draw=none] {\color{red} \ref*{thm:nat-succ}-\ref*{rule:all}}
edge from parent node[right, draw=none] {\color{red} 2}
                                  }
                         }
                    }
                                    edge from parent node[left, draw=none] {\color{red} \ref*{rule:exist}}
edge from parent node[right, draw=none] {\color{red} 0}
                  };
		\end{tikzpicture}
	\end{center}	

  \begin{theorem}[Upper bound for successor domain, \ref{def:id-ext-succ}, \ref{def:ind-nat}, \ref{def:successor-function}]\ \\ \label{thm:upper-successor}
$$ \lall N\xi \rall \lall \mem \eta {s(\xi)}\rall \lall \nonmem \eta \xi \rall \xi{\eq} \eta$$
\end{theorem}
	\begin{center}
		\begin{tikzpicture}[every node/.style = {draw, shape=rectangle,anchor=north},growth parent anchor=south,level distance=.4 cm,sibling distance = 4.7cm]
		\node{$ \lexi N\xi \rexi \lexi \mem \eta {s(\xi)}\rexi \lexi \nonmem \eta \xi \rexi \xi\ds{\eq} \eta$}
 			child{[level distance=.1cm]node{$ Nx $}
				child{[level distance=.1cm]node{$ \mem y {s(x)} $}
				child{[level distance=.4cm]node{$ \lexi \nonmem y x \rexi x\ds{\eq} y$}
                                 child{[level distance=.4cm]node{$ Hs(x)x $}
                                  child{[level distance=.4cm]node{$ \lall \mem \xi  {s(x)}\rall  Ex\xi$}
                                  child{[level distance=.4cm]node{$ Exy $}
                                  child{[level distance=.4cm]node{$ \lall \nonmem y x\rall x=y $}
                                  child{[level distance=.4cm]node{$  $}
edge from parent node[right, draw=none] {\color{red} 4}
                                  }
                                    edge from parent node[left, draw=none] {\color{red} \ref*{def:id-ext-succ}-\ref*{rule:subs}}
edge from parent node[right, draw=none] {\color{red} 0}
                                  }
                                    edge from parent node[left, draw=none] {\color{red} \ref*{rule:all}}
edge from parent node[right, draw=none] {\color{red} 0-3}
                                  }
                                   edge from parent node[left, draw=none] {\color{red} \ref*{def:id-ext-succ}-\ref*{rule:subs}-\ref*{rule:exist}}
edge from parent node[right, draw=none] {\color{red} 0}
                                  }
                                    edge from parent node[left, draw=none] {\color{red} \ref*{thm:nat-succ}}
edge from parent node[right, draw=none] {\color{red} 2}
                                }
                         }
                    }
                                    edge from parent node[left, draw=none] {\color{red} \ref*{rule:exist}}
edge from parent node[right, draw=none] {\color{red} 0}
                  };
		\end{tikzpicture}
	\end{center}	

We pause to observe by induction that zero in in the successor of every natural number.

  \begin{theorem}[Zero in successor, \ref{def:empty-object}, \ref{def:ind-nat}, \ref{def:successor-function}]\ \\ \label{thm:zero-in-suc}
  $$ \lall N\xi \rall  \mem o {s(\xi)} $$
\end{theorem}

\begin{center}
		\begin{tikzpicture}[every node/.style = {draw, shape=rectangle,anchor=north},growth parent anchor=south,level distance=.4 cm,sibling distance = 5cm]
		\node{$ \lexi N\xi \rexi \nonmem o {s(\xi)}$}
				child{[level distance=.1cm]node{$ A\xi \prop \mem o {s(\xi)}$}
                                 child{[level distance=.4cm]node{$  \lexi N\xi \rexi \lexi A\xi \rexi \ds{A}s(\xi) $}
                                  child{[level distance=.1cm]node{$ Nx $}
                                  child{[level distance=.1cm]node{$ Ax$}
                                  child{[level distance=.4cm]node{$ \ds{A}s(x)$}
                                  child{[level distance=.4cm]node{$ \mem o {s(x)}$}
                                  child{[level distance=.4cm]node{$ \nonmem o {s(s(x))}$}
                                  child{[level distance=.4cm]node{$ \mem o {s(s(x))}$}
                                  child{[level distance=.4cm]node{$ $}
edge from parent node[right, draw=none] {\color{red} 1}
                                  }
                                                                     edge from parent node[left, draw=none] {\color{red} \ref*{thm:growth-nat-num}}
edge from parent node[right, draw=none] {\color{red} 1-4}
                                  }
                                                                      edge from parent node[left, draw=none] {\color{red} \ref*{rule:subs}}
edge from parent node[right, draw=none] {\color{red} 1-5}
                                  }
                                    edge from parent node[left, draw=none] {\color{red} \ref*{rule:subs}}
edge from parent node[right, draw=none] {\color{red} 1-4}
                                  }
                                  }
                                  }
                                    edge from parent node[left, draw=none] {\color{red} \ref*{rule:exist}}
edge from parent node[right, draw=none] {\color{red} 0}
                                  }
                                  }
                                  child{[level distance=.1cm, sibling distance=3cm]node{$ \lall N\xi \rall \lall A\xi \rall As(\xi) $}
                                  child{[level distance=.4cm]node{$  \ds{A}o $}
                                  child{[level distance=.4cm]node{$ \nonmem o {s(o)} $}
                                  child{[level distance=.4cm]node{$ \mem o {s(o)} $}
                                  child{[level distance=.4cm]node{$ $}
edge from parent node[right, draw=none] {\color{red} 1}
                                  }
                                    edge from parent node[left, draw=none] {\color{red}  \ref*{thm:peanoi}-\ref*{thm:nat-num-in-suc}}
                                  }
                                    edge from parent node[left, draw=none] {\color{red} \ref*{rule:subs}}
edge from parent node[right, draw=none] {\color{red} 0-2}
                                  }
                                  }
                                  child{[level distance=.4cm]node{$ Ao $}
  			child{[level distance=.1cm]node{$ Nx $}
				child{[level distance=.4cm]node{$ \nonmem o {s(x)} $}
                                 child{[level distance=.4cm]node{$  A x $}
                                  child{[level distance=.4cm]node{$ \mem o {s(x)} $}
                                  child{[level distance=.4cm]node{$ $}
edge from parent node[right, draw=none] {\color{red} 2}
                                  }
                                  edge from parent node[left, draw=none] {\color{red} \ref*{rule:subs} }
edge from parent node[right, draw=none] {\color{red} 0-5}
                                  }
                                   edge from parent node[left, draw=none] {\color{red} \ref*{thm:peanov}-\ref*{rule:all}}
 edge from parent node[right, draw=none] {\color{red} 1-2-3}
                                  }
                                  }
                                                                     edge from parent node[left, draw=none] {\color{red} \ref*{rule:exist}}
edge from parent node[right, draw=none] {\color{red} 3}
                          }
                         }
                    }
                                    edge from parent node[left, draw=none] {\color{red} \ref*{rule:con-pro}}
                  };
		\end{tikzpicture}
	\end{center}	

The next theorem establishes that containment is transitive across the natural numbers.

  \begin{theorem}[Transitivity of containment, \ref{def:empty-object}, \ref{def:ind-nat}, \ref{def:successor-function}]\ \\ \label{thm:transitivity-containment}
    $$ \lall N\eta \rall \lall N\zeta \rall \lall \mem \xi \eta \rall 
    \lall \mem \eta \zeta\rall \mem \xi \zeta $$
    \end{theorem}
    \begin{center}
		\begin{tikzpicture}[every node/.style = {draw, shape=rectangle,anchor=north},growth parent anchor=south,level distance=.4 cm,sibling distance = 4.5cm]
                  \node{$ \lexi N\eta \rexi \lexi N\zeta \rexi \lexi \mem \xi \eta \rexi  \lexi \mem \eta \zeta\rexi \nonmem \xi \zeta $}
				child{[level distance=.1cm]node{$ Ny $}
				child{[level distance=.1cm]node{$ Nz $}
				child{[level distance=.1cm]node{$ \mem x y $}
                                  child{[level distance=.4cm]node{$ \lexi \mem y z\rexi \nonmem x z $}
                                  child{[level distance=.1cm]node{$ A\zeta \prop \lall \mem y \zeta \rall \mem x \zeta  $}
                                  child{[level distance=.4cm]node{$ \lexi N\xi \rexi \lexi A\xi \rexi \ds{A}s(\xi) $}
                                  child{[level distance=.1cm]node{$ Nu $}
                                  child{[level distance=.1cm]node{$ Au$}
                                  child{[level distance=.4cm]node{$ \ds{A}s(u) $}
                                  child{[level distance=.4cm]node{$ \lall \mem y u \rall \mem x u  $}
                                  child{[level distance=.1cm, sibling distance=3.1cm]node{$ \lexi \mem y {s(u)}  \rexi \nonmem x {s(u)}$}
                                  child{[level distance=.4cm]node{$ \mem y u$}
                                  child{[level distance=.4cm]node{$ \mem x u$}
                                  child{[level distance=.4cm]node{$ \mem x {s(u)}$}
                                  child{[level distance=.4cm]node{$ \nonmem x {s(u)}$}
                                  child{[level distance=.4cm]node{$ $}
edge from parent node[right, draw=none] {\color{red} 1}
                                  }
                                    edge from parent node[left, draw=none] {\color{red} \ref*{rule:exist}}
edge from parent node[right, draw=none] {\color{red} 3}
                                  }
                                    edge from parent node[left, draw=none] {\color{red} \ref*{thm:growth-nat-num}-\ref*{rule:all}}
edge from parent node[right, draw=none] {\color{red} 0-6}
                                  }
                                    edge from parent node[left, draw=none] {\color{red} \ref*{rule:all}}
edge from parent node[right, draw=none] {\color{red} 0-2}
                                  }
                                  }
                                  child{[level distance=.4cm]node{$ \nonmem y u$}
                                  child{[level distance=.1cm]node{$ \mem y {s(u)}$}
                                  child{[level distance=.4cm]node{$ \nonmem x {s(u)}$}
                                  child{[level distance=.4cm]node{$ y=u$}
                                  child{[level distance=.4cm]node{$ \nonmem x  {s(y)}$}
                                  child{[level distance=.4cm]node{$ \mem x {s(y)} $}
                                  child{[level distance=.4cm]node{$  $}
edge from parent node[right, draw=none] {\color{red} 1}
                                  }
                                    edge from parent node[left, draw=none] {\color{red} \ref*{thm:growth-nat-num}-\ref*{rule:all}}
edge from parent node[right, draw=none] {\color{red} 13-15}
                                  }
                                    edge from parent node[left, draw=none] {\color{red} \ref*{rule:subf}}
edge from parent node[right, draw=none] {\color{red} 0-1}
                                  }
                                    edge from parent node[left, draw=none] {\color{red} \ref*{thm:upper-successor}-\ref*{rule:all}}
edge from parent node[right, draw=none] {\color{red} 1-2-7}
                                  }
                                  }
                                    edge from parent node[left, draw=none] {\color{red} \ref*{rule:exist}}
edge from parent node[right, draw=none] {\color{red} 1}
                                  }
                                  }
                                    edge from parent node[left, draw=none] {\color{red} \ref*{rule:subs}}
edge from parent node[right, draw=none] {\color{red} 1-5}
                                  }
                                    edge from parent node[left, draw=none] {\color{red} \ref*{rule:subs}}
edge from parent node[right, draw=none] {\color{red} 1-4}
                                  }
                                  }
                                  }
                                    edge from parent node[left, draw=none] {\color{red} \ref*{rule:exist}}
edge from parent node[right, draw=none] {\color{red} 0}
                                  }
                                  }
                                  child{[level distance=.1cm, sibling distance=2.8cm]node{$  \lall N\xi \rall \lall A\xi \rall As(\xi)$}
                                  child{[level distance=.4cm]node{$ Ao $}
                                  child{[level distance=.4cm]node{$ Az $}
                                  child{[level distance=.4cm]node{$ \lall \mem y z \rall \mem x z $}
                                  child{[level distance=.4cm]node{$  $}
edge from parent node[right, draw=none] {\color{red} 5}
                                  }
                                    edge from parent node[left, draw=none] {\color{red} \ref*{rule:subs}}
edge from parent node[right, draw=none] {\color{red} 0-3}
                                  }
                                    edge from parent node[left, draw=none] {\color{red} \ref*{thm:peanov}-\ref*{rule:all}}
edge from parent node[right, draw=none] {\color{red} 0-1-5}
                                  }
                                  }
                                  child{[level distance=.4cm]node{$ \ds{A}o $}
                                  child{[level distance=.4cm]node{$ \lexi \mem y o \rexi \nonmem x o $}
                                  child{[level distance=.4cm]node{$ Oo $}
                                  child{[level distance=.4cm]node{$ \lall \xi=\xi \rall \nonmem \xi o$}
                                child{[level distance=.4cm]node{$ \nonmem y o $}
                                  child{[level distance=.4cm]node{$ \mem y o $}
                                  child{[level distance=.4cm]node{$ $}
edge from parent node[right, draw=none] {\color{red} 1}
                                  }
                                    edge from parent node[left, draw=none] {\color{red} \ref*{rule:exist}}
edge from parent node[right, draw=none] {\color{red} 3}
                                  }
                                    edge from parent node[left, draw=none] {\color{red} \ref*{rule:contradiction}-\ref*{rule:all}}
edge from parent node[right, draw=none] {\color{red} 0}
                                  }
                                  edge from parent node[left, draw=none] {\color{red} \ref*{def:empty}-\ref*{def:empty-object}-\ref*{rule:subs}}
edge from parent node[right, draw=none] {\color{red} 0}
                                  }
                                    edge from parent node[left, draw=none] {\color{red} \ref*{def:empty-object}}
                                  }
                                    edge from parent node[left, draw=none] {\color{red} \ref*{rule:subs}}
edge from parent node[right, draw=none] {\color{red} 0-2}
                                  }
                                  }
                                  }
                                    edge from parent node[left, draw=none] {\color{red} \ref*{rule:con-pro}}
                                }
                         }
                    }
                    }
                                    edge from parent node[left, draw=none] {\color{red} \ref*{rule:exist}}
edge from parent node[right, draw=none] {\color{red} 0}
                  };
		\end{tikzpicture}
	\end{center}	

Finally, we are in a position to prove the third Peano axiom. It states uniqueness for the predecessor of any natural number.

  \begin{theorem}[Peano axiom III, \ref{def:ind-nat}, \ref{def:successor-function}]\ \\ \label{thm:peanoiii}
  $$ \lall N\xi \rall \lall N\eta \rall \lall s(\xi)=s(\eta) \rall \xi=\eta $$
\end{theorem}

\begin{center}
		\begin{tikzpicture}[every node/.style = {draw, shape=rectangle,anchor=north},growth parent anchor=south,level distance=.4 cm,sibling distance = 5cm]
		\node{$ \lexi N\xi \rexi \lexi N\eta \rexi \lexi s(\xi)=s(\eta) \rexi \xi\ds{\eq} \eta$}
 			child{[level distance=.1cm]node{$ Nx $}
				child{[level distance=.1cm]node{$ Ny $}
				child{[level distance=.1cm]node{$ s(x)=s(y)$}
                                  child{[level distance=.1cm]node{$ x\ds{\eq} y $}
                                  child{[level distance=.4cm]node{$  \nonmem y x $}
                                  child{[level distance=.4cm]node{$ \mem y {s(y)} $}
                                  child{[level distance=.4cm]node{$ \mem y {s(x)}$}
                                  child{[level distance=.4cm]node{$ y=x $}
                                  child{[level distance=.1cm]node{$ x\ds{\eq} x$}
                                  child{[level distance=.4cm]node{$ $}
                                  }
                                    edge from parent node[left, draw=none] {\color{red} \ref*{rule:subf}}
edge from parent node[right, draw=none] {\color{red} 0-4}
                                  }
                                    edge from parent node[left, draw=none] {\color{red} \ref*{thm:upper-successor}-\ref*{rule:all}}
edge from parent node[right, draw=none] {\color{red} 0-2-6}
                                  }
                                    edge from parent node[left, draw=none] {\color{red}  \ref*{rule:subf}}
edge from parent node[right, draw=none] {\color{red} 0-3}
                                  }
                                    edge from parent node[left, draw=none] {\color{red} \ref*{thm:nat-num-in-suc}-\ref*{rule:all}}
edge from parent node[right, draw=none] {\color{red} 3}
                                  }
                                  }
                                  child{[level distance=.1cm, sibling distance=3cm]node{$  \mem y x $}
                                  child{[level distance=.4cm]node{$  \nonmem x y $}
                                  child{[level distance=.4cm]node{$ \mem x {s(x)} $}
                                  child{[level distance=.4cm]node{$ \mem x {s(y)}$}
                                  child{[level distance=.4cm]node{$ x=y $}
                                  child{[level distance=.4cm]node{$ $}
edge from parent node[right, draw=none] {\color{red} 5}
                                  }
                                    edge from parent node[left, draw=none] {\color{red} \ref*{thm:upper-successor}-\ref*{rule:all}}
edge from parent node[right, draw=none] {\color{red} 0-2-6}
                                  }
                                    edge from parent node[left, draw=none] {\color{red}  \ref*{rule:subf}}
edge from parent node[right, draw=none] {\color{red} 0-4}
                                  }
                                    edge from parent node[left, draw=none] {\color{red} \ref*{thm:nat-num-in-suc}-\ref*{rule:all}}
edge from parent node[right, draw=none] {\color{red} 5}
                                  }
                                  }
                                  child{[level distance=.4cm]node{$ \mem x y $}
                                  child{[level distance=.4cm]node{$  \mem x x $}
                                  child{[level distance=.4cm]node{$ \nonmem x x $}
                                  child{[level distance=.4cm]node{$ $}
edge from parent node[right, draw=none] {\color{red} 1}
                                  }
                                  edge from parent node[left, draw=none] {\color{red} \ref*{thm:irreflexive-containment}-\ref*{rule:all}}
edge from parent node[right, draw=none] {\color{red} 6}
                                  }
                                   edge from parent node[left, draw=none] {\color{red} \ref*{thm:transitivity-containment}-\ref*{rule:all}}
 edge from parent node[right, draw=none] {\color{red} 0-1-4-5}
                                  }
                                  }
                                  }
                          }
                         }
                    }
                                    edge from parent node[left, draw=none] {\color{red} \ref*{rule:exist}}
edge from parent node[right, draw=none] {\color{red} 0}
                  };
		\end{tikzpicture}
	\end{center}	

\end{document}